\documentclass[10pt]{amsart}
\usepackage{amsfonts}
\usepackage{amsmath}
\usepackage{amsthm}
\usepackage{amssymb}
\usepackage{latexsym}
\newtheorem{cor}{Corollary}[section]
\newtheorem{lem}{Lemma}[section]

\newtheorem{prop}{Proposition}[section]

\theoremstyle{definition}
\newtheorem{defn}{Definition}[section]
\theoremstyle{definition}
\newtheorem{thm}{Theorem}

\newenvironment{pf}{\proof}{\endproof}

\theoremstyle{remark}
\newtheorem{remark}{Remark}[section]

\numberwithin{equation}{section}
\begin{document}

\newcommand{\thmref}[1]{Theorem~\ref{#1}}
\newcommand{\secref}[1]{Sect.~\ref{#1}}
\newcommand{\lemref}[1]{Lemma~\ref{#1}}
\newcommand{\propref}[1]{Proposition~\ref{#1}}
\newcommand{\corref}[1]{Corollary~\ref{#1}}
\newcommand{\remref}[1]{Remark~\ref{#1}}
\newcommand{\nc}{\newcommand}
\newcommand{\rnc}{\renewcommand}
\nc{\cal}{\mathcal}
\nc{\goth}{\mathfrak}
\rnc{\bold}{\mathbf}
\renewcommand{\frak}{\mathfrak}
\renewcommand{\Bbb}{\mathbb}

\nc{\Cal}{\mathcal}
\nc{\Xp}[1]{X^+(#1)}
\nc{\Xm}[1]{X^-(#1)}
\nc{\on}{\operatorname}
\nc{\ch}{\mbox{ch}}
\nc{\Z}{{\bold Z}}
\nc{\J}{{\cal J}}
\nc{\C}{{\bold C}}
\nc{\Q}{{\bold Q}}
\renewcommand{\P}{{\cal P}}
\nc{\N}{{\Bbb N}}
\nc\beq{\begin{equation}}
\nc\enq{\end{equation}}
\nc\lan{\langle}
\nc\ran{\rangle}
\nc\bsl{\backslash}
\nc\mto{\mapsto}
\nc\lra{\leftrightarrow}
\nc\hra{\hookrightarrow}
\nc\sm{\smallmatrix}
\nc\esm{\endsmallmatrix}
\nc\sub{\subset}
\nc\ti{\tilde}
\nc\nl{\newline}
\nc\fra{\frac}
\nc\und{\underline}
\nc\ov{\overline}
\nc\ot{\otimes}
\nc\bbq{\bar{\bq}_l}
\nc\bcc{\thickfracwithdelims[]\thickness0}
\nc\ad{\text{\rm ad}}
\nc\Ad{\text{\rm Ad}}
\nc\Hom{\text{\rm Hom}}
\nc\End{\text{\rm End}}
\nc\Ind{\text{\rm Ind}}
\nc\Res{\text{\rm Res}}
\nc\Ker{\text{\rm Ker}}
\rnc\Im{\text{Im}}
\nc\sgn{\text{\rm sgn}}
\nc\tr{\text{\rm tr}}
\nc\Tr{\text{\rm Tr}}
\nc\supp{\text{\rm supp}}
\nc\card{\text{\rm card}}
\nc\bst{{}^\bigstar\!}
\nc\he{\heartsuit}
\nc\clu{\clubsuit}
\nc\spa{\spadesuit}
\nc\di{\diamond}

\nc\al{\alpha}
\nc\bet{\beta}
\nc\ga{\gamma}
\nc\de{\delta}
\nc\ep{\epsilon}
\nc\io{\iota}
\nc\om{\omega}
\nc\si{\sigma}
\rnc\th{\theta}
\nc\ka{\kappa}
\nc\la{\lambda}
\nc\ze{\zeta}

\nc\vp{\varpi}
\nc\vt{\vartheta}
\nc\vr{\varrho}

\nc\Ga{\Gamma}
\nc\De{\Delta}
\nc\Om{\Omega}
\nc\Si{\Sigma}
\nc\Th{\Theta}
\nc\La{\Lambda}
\nc\boa{\bold a}
\nc\bob{\bold b}
\nc\boc{\bold c}
\nc\bod{\bold d}
\nc\boe{\bold e}
\nc\bof{\bold f}
\nc\bog{\bold g}
\nc\boh{\bold h}
\nc\boi{\bold i}
\nc\boj{\bold j}
\nc\bok{\bold k}
\nc\bol{\bold l}
\nc\bom{\bold m}
\nc\bon{\bold n}
\nc\boo{\bold o}
\nc\bop{\bold p}
\nc\boq{\bold q}
\nc\bor{\bold r}
\nc\bos{\bold s}
\nc\bou{\bold u}
\nc\bov{\bold v}
\nc\bow{\bold w}
\nc\boz{\bold z}

\nc\ba{\bold A}
\nc\bb{\bold B}
\nc\bc{\bold C}
\nc\bd{\bold D}
\nc\be{\bold E}
\nc\bg{\bold G}
\nc\bh{\bold h}
\nc\bH{\bold H}

\nc\bi{\bold I}
\nc\bj{\bold J}
\nc\bk{\bold K}
\nc\bl{\bold L}
\nc\bm{\bold M}
\nc\bn{\bold N}
\nc\bo{\bold O}
\nc\bp{\bold P}
\nc\bq{\bold Q}
\nc\br{\bold R}
\nc\bs{\bold S}
\nc\bt{\bold T}
\nc\bu{\bold U}
\nc\bv{\bold v}
\nc\bV{\bold V}

\nc\bw{\bold W}
\nc\bz{\bold Z}
\nc\bx{\bold x}
\nc\bX{\bold X}
\nc\blambda{{\mbox{\boldmath $\Lambda$}}}
\nc\bpi{{\mbox{\boldmath $\pi$}}}

\nc\e[1]{E_{#1}}
\nc\ei[1]{E_{\delta - \alpha_{#1}}}
\nc\esi[1]{E_{s \delta - \alpha_{#1}}}
\nc\eri[1]{E_{r \delta - \alpha_{#1}}}
\nc\ed[2][]{E_{#1 \delta,#2}}
\nc\ekd[1]{E_{k \delta,#1}}
\nc\emd[1]{E_{m \delta,#1}}
\nc\erd[1]{E_{r \delta,#1}}

\nc\ef[1]{F_{#1}}
\nc\efi[1]{F_{\delta - \alpha_{#1}}}
\nc\efsi[1]{F_{s \delta - \alpha_{#1}}}
\nc\efri[1]{F_{r \delta - \alpha_{#1}}}
\nc\efd[2][]{F_{#1 \delta,#2}}
\nc\efkd[1]{F_{k \delta,#1}}
\nc\efmd[1]{F_{m \delta,#1}}
\nc\efrd[1]{F_{r \delta,#1}}
\nc{\ug}{\bu^{fin}}

\nc\fa{\frak a}
\nc\fb{\frak b}
\nc\fc{\frak c}
\nc\fd{\frak d}
\nc\fe{\frak e}
\nc\ff{\frak f}
\nc\fg{\frak g}
\nc\fh{\frak h}
\nc\fj{\frak j}
\nc\fk{\frak k}
\nc\fl{\frak l}
\nc\fm{\frak m}
\nc\fn{\frak n}
\nc\fo{\frak o}
\nc\fp{\frak p}
\nc\fq{\frak q}
\nc\fr{\frak r}
\nc\fs{\frak s}
\nc\ft{\frak t}
\nc\fu{\frak u}
\nc\fv{\frak v}
\nc\fz{\frak z}
\nc\fx{\frak x}
\nc\fy{\frak y}

\nc\fA{\frak A}
\nc\fB{\frak B}
\nc\fC{\frak C}
\nc\fD{\frak D}
\nc\fE{\frak E}
\nc\fF{\frak F}
\nc\fG{\frak G}
\nc\fH{\frak H}
\nc\fJ{\frak J}
\nc\fK{\frak K}
\nc\fL{\frak L}
\nc\fM{\frak M}
\nc\fN{\frak N}
\nc\fO{\frak O}
\nc\fP{\frak P}
\nc\fQ{\frak Q}
\nc\fR{\frak R}
\nc\fS{\frak S}
\nc\fT{\frak T}
\nc\fU{\frak U}
\nc\fV{\frak V}
\nc\fZ{\frak Z}
\nc\fX{\frak X}
\nc\fY{\frak Y}
\nc\tfi{\ti{\Phi}}
\nc\bF{\bold F}

\nc\ua{\bold U_\A}

\nc\qinti[1]{[#1]_i}
\nc\q[1]{[#1]_q}
\nc\xpm[2]{E_{#2 \delta \pm \alpha_#1}}  
\nc\xmp[2]{E_{#2 \delta \mp \alpha_#1}}
\nc\xp[2]{E_{#2 \delta + \alpha_{#1}}}
\nc\xm[2]{E_{#2 \delta - \alpha_{#1}}}
\nc\hik{\ed{k}{i}}
\nc\hjl{\ed{l}{j}}
\nc\qcoeff[3]{\left[ \begin{smallmatrix} {#1}& \\ {#2}& \end{smallmatrix}
\negthickspace \right]_{#3}}
\nc\qi{q}
\nc\qj{q}

\nc\ufdm{{_\ca\bu}_{\rm fd}^{\le 0}}


\nc\isom{\cong} 

\nc{\pone}{{\Bbb C}{\Bbb P}^1}
\nc{\pa}{\partial}
\def\H{\cal H}
\def\L{\cal L}
\nc{\F}{{\cal F}}
\nc{\Sym}{{\goth S}}
\nc{\A}{{\cal A}}
\nc{\arr}{\rightarrow}
\nc{\larr}{\longrightarrow}

\nc{\ri}{\rangle}
\nc{\lef}{\langle}
\nc{\W}{{\cal W}}
\nc{\uqatwoatone}{{U_{q,1}}(\su)}
\nc{\uqtwo}{U_q(\goth{sl}_2)}
\nc{\dij}{\delta_{ij}}
\nc{\divei}{E_{\alpha_i}^{(n)}}
\nc{\divfi}{F_{\alpha_i}^{(n)}}
\nc{\Lzero}{\Lambda_0}
\nc{\Lone}{\Lambda_1}
\nc{\ve}{\varepsilon}
\nc{\phioneminusi}{\Phi^{(1-i,i)}}
\nc{\phioneminusistar}{\Phi^{* (1-i,i)}}
\nc{\phii}{\Phi^{(i,1-i)}}
\nc{\Li}{\Lambda_i}
\nc{\Loneminusi}{\Lambda_{1-i}}
\nc{\vtimesz}{v_\ve \otimes z^m}

\nc{\asltwo}{\widehat{\goth{sl}_2}}
\nc\eh{\frak h^e}  
\nc\loopg{L(\frak g)}  
\nc\eloopg{L^e(\frak g)} 
\nc\ebu{\bu^e} 

\nc\teb{\tilde E_\boc}
\nc\tebp{\tilde E_{\boc'}}

\title{Weyl modules for classical and quantum affine algebras }
\author{Vyjayanthi Chari}
\address{Vyjayanthi Chari, Department of Mathematics, University of California, 
Riverside, CA 92521.}
\author{Andrew Pressley}
\address{Andrew Pressley, Department of Mathematics, Kings College, London, WC 
2R, 2LS, England, U.K.}

\pagestyle{plain} \maketitle
\section{Introduction} The study of the irreducible finite--dimensional 
representations of quantum affine algebras has been the subject of a number of 
papers, \cite{AK}, \cite{CP3}, \cite{CP5}, \cite{FR}, \cite{FM}, \cite{GV}, 
\cite{KS} to name a few. However, the structure of these representations is 
still unknown except in certain special cases. In this paper, we approach the 
problem by  studying the classical ($q\to 1$) limits of these representations. 
Standard results imply, for example, that if $V$ is  a finite--dimensional 
representation of $\bu_q(\hat{\frak g})$, its $q\to 1$ limit $\overline V$ has 
the same structure  as  a $\frak g$--module as $V$ has as a  $U_q(\frak 
g)$--module.

We begin by studying an appropriate class of representations of the affine Lie 
algebra $\hat{\frak g}$. The finite--dimensional irreducible representations of 
$\hat{\frak g}$ were classified in \cite{C}, \cite{CP1}, where it was shown that 
such representations are highest weight in a suitable sense, the highest weight 
being an $n$--tuple of polynomials $\bpi$, where $n$ is the rank of $\frak g$.
 We therefore study the class of all highest weight finite--dimensional 
representations of $\hat{\frak g}$. In fact, we prove that corresponding to each 
irreducible finite--dimensional representation $V(\bpi)$ there exists a unique  
(up to isomorphism) finite--dimensional highest weight module  $W(\bpi)$, such 
that {\it any} finite dimensional highest weight module $V$ with highest weight 
$\bpi$ is a quotient of $W(\bpi)$. We call these modules the Weyl modules 
because of an analogy with the modular representation theory of $\frak g$, which we now explain.

In \cite{CP5}, we showed  that the  irreducible representations of 
$U_q(\hat{\frak g})$ are also highest weight and that their isomorphism classes 
are parametrized by a $n$--tuples of polynomials $\bpi_q$ with coefficients in 
$\bc(q)$. Under a natural condition on $\bpi_q$, the corresponding 
representation $V_q(\bpi_q)$ of $\bu_q(\hat{\frak g})$ specializes as $q\to 1$ 
to a representation $\overline{V_q(\bpi_q)} $ of $\hat{\frak g}$ and is a 
quotient of $W(\bpi)$, where $\bpi$ is obtained from $\bpi_q$ by setting $q=1$. 
We conjecture that every $W(\bpi)$ is the classical limit of an irreducible  
$U_q(\hat{\frak g})$-module. This is analogous to the fact that the Weyl modules 
for $\frak g$ in characteristic $p$ are the mod $p$ reductions of the 
irreducible modules in characteristic zero. We prove the conjecture  in the case 
of $\frak 
g=sl_2$ in this paper. The conjecture is also true for the fundamental 
representations of $\bu_q(\hat{\frak g})$, but the details of that will appear 
elsewhere.  

In Section 3, we prove a factorization property of Weyl modules analogous the 
one for the irreducible modules proved in \cite{CP1}. In Section 5We obtain a necessary 
and sufficient condition for the Weyl modules to be irreducible: the interesting feature of this proof is that it uses the fact that the specialized irreducible 
modules for the quantum algebra are quotients of the Weyl module.  Further, the 
condition for irreducibility of the Weyl modules is the same as a condition that appears first in the work of Drinfeld on the closely related Yangians, 
\cite{Dr1}.

The Weyl modules we define are quotients of a family of level zero integrable 
modules for the extended affine Lie algebra, one corresponding to each dominant 
integral weight of $\frak g$. We call these modules $W(\lambda)$. According to 
unpublished work of Kashiwara, these modules are the classical analogues of the 
modules $V^{max}(\lambda)$ defined in \cite{K}. Further, Kashiwara has a number 
of conjectures on the crystal basis of $V^{max}(\lambda)$.    In Section 6, we 
identify the modules $W(\lambda)$ explicitly  in the case of $sl_2$. A similar identification can then be proved for the modules $V^{max}(\lambda)$  which settles one of Kashiwara's conjectures, but the details of this will appear elsewhere. 

\vskip9pt\noindent{\it{Acknowledgements.}} We thank F.Knop, D.Rush and R.Sujatha for discussions. We also thank  M.Kashiwara for explaining his conjectures to us and for drawing our attention to the modules $V^{max}(\lambda)$.

\section {Preliminaries and Some Identities} Let  $\frak{g}$ be a 
finite-dimensional complex simple Lie algebra of rank $n$,  $\frak h$  a Cartan 
subalgebra of $\frak g$ and  $R$  the set of roots of $\frak g$ with 
respect to $\frak h$. Let $I=\{1,2,\cdots ,n\}$, fix a set of simple 
roots $\alpha_i$ ($i\in I$), and let $R^+$ be the corresponding set of positive 
roots.  Let $\theta\in R^+$ be the highest root in $R^+$.  For 
$\alpha\in R^+$, fix non-zero elements  $x_\alpha^\pm\in\frak g$, 
$h_\alpha\in\frak h$ such that 
\begin{equation*}
[x_\alpha^+, x_\alpha^-] =h_\alpha,\ \  [h_\alpha,x_\alpha^\pm] = \pm 2 
x_\alpha^\pm .
\end{equation*}  
Let $Q=\bigoplus_{i=1}^n \bz\alpha_i$ (resp. 
$Q_+=\bigoplus_{i=1}^n \bn\alpha_i$) denote the root (resp. positive root) 
lattice 
of $\frak g$.  For $\eta\in Q^+$, $\eta=\sum_ir_i\alpha_i$, we set 
${\text{ht}}\,\eta=\sum_i r_i$.
 The lattice $P$ (resp. $P_+$) of integral (resp. dominant 
integral) weights is the set of elements $\lambda\in\frak h^*$ such that 
$\lambda(h_\alpha)\in\bz$ for all $\alpha\in R$ (resp. 
$\lambda(h_\alpha)\ge 0$ for all $\alpha\in R^+$). For $i\in I$,  the 
fundamental weight 
$\omega_i$ of $\frak g$ is given by $\omega_i(\alpha_j)=\delta_{ij}$.  Let 
$<\ ,\ >$ be the bilinear pairing on $P$ such that $<\alpha_i,\omega_j> 
=\delta_{ij}$. Set 
$a_{ij}=<\alpha_i,\alpha_j>$ for $i,j\in I$. The bilinear form  induces an 
isomorphism $\frak h\cong \frak h^*$ such that, 
if $\beta=\sum_i r_i\alpha_i\in R^+$, then 
\begin{equation*} 
h_\beta = \sum_j \frac{d_jr_j}{d_\beta}h_j, 
\end{equation*}
where for a root $\alpha\in R$, we set $d_\alpha =\frac12<\alpha,\alpha>$.   
Let  
$W\subset {\text{Aut}}(\frak h^*)$ be the Weyl group of $\frak g$; it is 
well known that $W$ is generated by simple reflections $s_i$ ($i\in I$).

The  extended loop algebra of $\frak g$ is the Lie algebra
\begin{equation*}\eloopg = \frak g\otimes \bc[t,t^{-1}]\oplus \bc 
d,\end{equation*}
with commutator given by
\begin{equation*}
[d, x\otimes t^r] = rx\otimes t^r,  \ \ \ [x\otimes t^r, y\otimes t^s] 
=[x,y]\otimes t^{r+s}\end{equation*}
for $x,y\in\frak g$, $r,s\in\bz$. The loop algebra $\loopg$ is the subalgebra 
$\frak g\otimes \bc[t,t^{-1}]$ of 
$\eloopg$. 
Let $\eh =\frak h\oplus\bc d$. Define $\delta\in(\eh)^*$ by
\begin{equation*}
\delta(\frak h) =0,\ \ \delta(d)=1.
\end{equation*}
Extend $\lambda\in\frak h^*$ to an element of $(\eh)^*$ by setting 
$\lambda(d)=0$. Set $P^e = \bigoplus_{i=1}^n \bz\omega_i\oplus\bz\delta$, 
and define $P_+^e$ in the obvious way. We regard $W$ as acting on $(\eh)^*$ by 
setting $w(\delta)=\delta$ for all $w\in W$.


For any $x\in \frak g$, $m\in\bz$, we denote by $x_m$ the element $x\otimes 
t^m\in\eloopg$. Set $e_i^\pm =x_{\alpha_i}^\pm\otimes 1$ and $e_0^\pm 
=x_\theta^\mp\otimes t^{\pm 1}$. Then, the elements $e_i^\pm$ ($i=0,\cdots ,n$) 
and $d$ generate $\eloopg$.

For any Lie algebra $\frak a$, the universal enveloping algebra of $\frak a$ is 
denoted by $\bu(\frak a)$.   We set 
\begin{equation*}
\bu(\eloopg) =\ebu, \ \ \bu(\loopg)= \bu,\ \ \bu(\frak g)=\ug. 
\end{equation*}
Let $\bu(<)$ (resp. $\bu(>)$) be the subalgebra of $\bu$ generated by  
the $x^-_{\alpha_i, m}$, (resp. $x^+_{\alpha_i, m}$)  for 
$i\in I$, $m\in \bz$. Clearly, $x^-_{\alpha,m}\in\bu(<)$ (resp. 
$x^+_{\alpha,m}\in\bu(>)$) for all 
$\alpha\in R^+$, $m\in\bz$.  Set $\ug(<) =\bu(<)\cap\ug$ and define $\ug(>)$ 
similarly. Finally, let $\bu(0)$ be the subalgebra of $\bu$ generated by 
$h_{\alpha, m}$ for $\alpha\in R$, $m\in\bz$, $m\ne 0$. 
We have 
\begin{align*}
\ug &=\ug(<)\bu(\frak h)\ug(>),\\
\ebu &=\bu(<)\bu(0)\bu(\eh)\bu(>).
\end{align*} 

\begin{lem} The assignment
$T(x_{\alpha_i,m}^\pm) =x_{\alpha_i, m\pm 1}^\pm$, for $i\in I$, $m\in\bz$, 
defines an algebra automorphism of $\bu$.\hfill\qedsymbol
\end{lem}

We next recall some identities in $\ebu$, which are most conveniently stated 
using generating series.  Thus, for any  $\beta\in R^+$,  we 
introduce the following power series in an indeterminate $u$:
\begin{align*} 
\tilde{X}_\beta^-(u)=\sum_{m=-\infty}^\infty x_{\beta,m}^-u^{m+1},\ \ \  &\ \ \ 
X_\beta^-(u)=\sum_{m=1}^\infty x_{\beta ,m}^-u^m,\\
X_\beta^+(u)=\sum_{m=0}^\infty x_{\beta,m}^+u^m,\ \ \ &\ \ \ 
X_{\beta,0}^-(u) =\sum_{m=0}^\infty x_{\beta,m}^-u^{m+1},\\
\tilde{H}_\beta(u)=\sum_{m=-\infty}^\infty h_{\beta,m}u^{m+1},\ \ \ & \ \ \ 
\Lambda_\beta^\pm(u)=\sum_{m=0}^\infty \Lambda_{\beta,\pm m}u^m 
={\text{exp}}\left(-\sum_{k=1}^\infty\frac{h_{\beta,\pm k}}{k}u^k\right).
\end{align*}
Set $x^\pm_{\alpha_i}=x_i^\pm$, $x_{\alpha_i,m}^\pm=x_{i,m}^\pm$, and define 
$h_i$, $\Lambda_i^\pm$, etc., similarly.

The next lemma follows easily from the  definition of the $\Lambda_{i,m}$ 
($=\Lambda_{\alpha_i,m}$).

\begin{lem} The subalgebra $\bu(0)$ of $\bu$ is generated by the elements 
 $\Lambda_{i,m}$, for $i\in I$, $m\in\bz$. 
\hfill\qedsymbol
\end{lem} 

For any power series $f$ in $u$ with coefficients in an algebra $A$,  let $f_m$ 
be the  
coefficient of $u^m$ ($m\in\bz$) and let $f'$ denote the derivative of $f$ with 
respect to $u$. For $x\in\bu$, $r\in \bz^+$, set
\begin{equation*} x^{(r)} =\frac {x^r}{r!}.\end{equation*}
For an algebra $A$, let $A_+$ denote the augmentation ideal. The next result is 
a reformulation of a result of Garland, \cite{G}.
\begin{lem}{\label{gar}} Let $s\ge r\ge 1$, $\beta\in R^+$. 
\begin{enumerate}
\item[(i)] 
\begin{equation*}(x_{\beta,0}^+)^{(r)} (x_{\beta,1}^-)^{(s)} = 
(-1)^r(X_\beta^-(u)^{(s-r)}\Lambda_\beta^+(u))_s\mod \bu\bu(>)_+.
\end{equation*}
\item[(ii)] 
\begin{equation*}(x_{\beta,1}^+)^{(r)}(x_{\beta,0}^-)^{(s)} = 
(-1)^r(X_{\beta,0}^-(u)^{(s-r)}\Lambda_\beta^+(u))_s\mod \bu\bu(>)_+.
\end{equation*}
\end{enumerate}
\end{lem}
\begin{pf} We use the identity in \cite[Lemma 7.1]{G}, in the form given in 
\cite[Lemma 5.1]{CP6} for its quantum version, namely
\begin{align*}
&(x_{\beta,0}^+)^{(r)} (x_{\beta,1}^-)^{(s)}\\&=
\sum_{t=0}^r\sum_{m=0}^t\sum_{k=0}^m(-1)^t
\left(u^{-s+t}X_\beta^-(u)^{(s-t)}\right)_{t-m}\Lambda_{\beta,k}
\left(X_\beta^+(u)^{(r-t)}\right)_{m-k}.
\end{align*}
In the sum on the right-hand side of the equality, we get an element of 
$\bu\bu(>)_+$ unless $t=r$ and $m=k$, and so we have
\begin{align*} (x_{\beta,0}^+)^{(r)} (x_{\beta,1}^-)^{(s)}&=(-1)^r\sum_{m=0}^r
\left(u^{-s+r}X_\beta^-(u)^{(s-r)}\right)_{r-m}
\Lambda_{\beta,m}\mod\bu\bu(>)_+\\
&=\sum_{r=0}^m(-1)^r\left(X_\beta^-(u)^{(s-r)}\right)_{s-m}\Lambda^+_\beta(u)
_m 
\mod\bu\bu(>)_+\\
&=(-1)^r\left(X_\beta^-(u)^{(s-r)}\Lambda^+_\beta(u)\right)_s\mod\bu\bu(>)_+.
\end{align*}
The identity in (ii) follows from (i) by applying the automorphism $T:\bu\to 
\bu$.
\end{pf}

We conclude this section with some elementary properties of integrable 
$\eloopg$-modules.

A representation $V$ of $\ebu$ is called integrable if the Chevalley generators 
$e_i^\pm$, for $i= 0,1,\cdots ,n$, act locally nilpotently on $V$ and 
\begin{equation*} 
V=\bigoplus_{\lambda\in(\eh)^*} V_\lambda,
\end{equation*}
where 
\begin{equation*}
V_{\lambda}= \{v\in V: h.v=\lambda(h)v\ \ \forall \ \ h\in\eh \}.
\end{equation*}
It is well known that this implies that the elements $x^\pm_{\beta,m}$ act 
locally nilpotently on $V$ for all $\beta\in R^+$ and $m\in\bz$. Set
\begin{equation*} V_\lambda^+ =\{v\in V_\lambda:x_{\beta,m}^+.v =0\ \forall \ 
\beta\in R^+, m\in\bz\}.\end{equation*}
It is easy to see that, if $V$ is integrable, then
$V_\lambda\ne 0$ (resp. $V_\lambda^+\ne 0$) only if $\lambda\in P^e$ (resp. 
$\lambda\in P_+^e$). Further, if $v\in V_\lambda^+$, then 
\begin{align}
{\label{xloc}}(x^-_{\beta, m})^{\lambda(h_\beta)+1}.v &=0,\\
V_\lambda\ne 0\implies V_{w\lambda}\ne 0\ & \ \forall\ w\in W.
\end{align}

If $\lambda\in P_+$, let $V^{fin}(\lambda)$ be the finite-dimensional 
irreducible $\bu^{fin}$-module with highest weight $\lambda$. If $\lambda\in 
P^e_+$, the restriction of $\lambda$ to $\frak h^*$, also denoted by $\lambda$, 
is in $P_+$. If $V$ is an integrable $\bu^e$-module and $0\ne v\in V_\lambda^+$, 
then $\bu^{fin}.v$ is a $\bu^{fin}$-submodule of $V$ isomorphic to 
$V^{fin}(\lambda)$.

\begin{prop}{\label {1.2}} Let $V$ be an integrable $\ebu$-module.
Let $\lambda\in P_+^e$, let $0\ne v\in V_\lambda^+$ and let $\beta\in R^+$. 
Then:
\begin{enumerate} 
\item[(i)] $\Lambda_{\beta,m}.v =0$ for $m>\lambda(h_\beta)$;
\item[(ii)] for $r\ge 1$, $s>\lambda(h_\beta)$, 
\begin{equation*}
\left(X_\beta^-(u)^r\Lambda_\beta^+(u)\right)_s.v=0,\ \ \ \ 
\left(X_{\beta,0}^-(u)^r\Lambda_\beta^+(u)\right)_s.v=0;
\end{equation*}
\item[(iii)] for all $s\in\bz$,
\begin{equation*} 
\left(\tilde X^-_\beta(u)\Lambda_\beta^+(u)\right)_s.v=0,\ \ \ \ 
\left(\tilde H_\beta(u)\Lambda_\beta^+(u)\right)_s.v=0;
\end{equation*}
\item[(iv)] $\Lambda_{\beta, -m}.v =0$ for all $m>\lambda(h_\beta)$;
\item[(v)] for $0\le m\le \lambda(h_\beta)$,
\begin{equation*}\Lambda_{\beta,\lambda(h_\beta)}\Lambda_{\beta, -m}.v 
=\Lambda_{\beta,\lambda(h_\beta)-m}.v.\end{equation*}\end{enumerate}
\end{prop}
\begin{pf} 
(i) This follows by taking $r=s>\lambda(h_\beta)$ in Lemma \ref{gar}(i) 
and using equation (\ref{xloc}).

\noindent (ii) This  follows from Lemma \ref{gar} by replacing $r$ by $s-r$
and using equation (\ref{xloc}).

\noindent (iii) Taking $r=\lambda(h_\beta)$, $s=\lambda(h_\beta)+1$ in Lemma 
\ref{gar}(i) gives
\begin{equation}\label{2} \sum_{m=0}^{\lambda(h_\beta)} 
x_{\beta,m+1}^-\Lambda_{\beta,\lambda(h_\beta)-m}. v =0.\end{equation}
Applying $h_{\beta,k}$, for any $k\in\bz$, to the above equation and noting that 
$h_{\beta,k}.v\in 
V_{\lambda}^+$, we get 
\begin{equation*} \sum_{m=0}^{\lambda(h_\beta)} 
x_{\beta,k+m+1}^-\Lambda_{\beta,\lambda(h_\beta)-m}. v =0,\end{equation*}
which can be written as
\begin{equation*}\left(\tilde{X}_\beta^-(u)\Lambda_\beta^+(u)\right)_{k+1}.v 
=0.\end{equation*}

Applying $x_{\beta,-s-1}^+$, for $s\in\bz$,  to both sides of equation (\ref{2}) 
gives
\begin{equation*}
\sum_{m=0}^{\lambda(\beta)}h_{\beta,m-s}\Lambda_{\beta,\lambda(
h_\beta)-m}.v =0,
\end{equation*}
i.e.,
\begin{equation*} 
\sum_{m=-s}^{\lambda(h_\beta)-s}h_{\beta,m} 
\Lambda_{\beta,\lambda(h_\beta)-s-m}.v =0.
\end{equation*}
Replacing $s$ by $\lambda(h_\beta)-s+1$ and using part (i) of the lemma, one 
sees that this identity is equivalent to the second identity in (iii).

\noindent  (iv) and (v). During the remainder of this proof, write 
$\Lambda_\beta=\Lambda_\beta^+$,
$\tilde\Lambda_\beta(u) =\Lambda_\beta^-(u^{-1})$, so that
\begin{equation*} \tilde\Lambda_\beta(u) ={\text{exp}}\left(-\sum_{k=1}^\infty 
\frac{h_{\beta,-k}}{k}u^{-k}\right).\end{equation*}
Note that, as operators on $V_\lambda^+$,
\begin{equation*}\left(\lambda(h_\beta)-u\frac{\Lambda_\beta'}{\Lambda_\beta} 
+u\frac{\tilde\Lambda_\beta'}{\tilde\Lambda_\beta}\right)=\tilde{H}_\beta.
\end{equation*}
By (iii), we have 
\begin{equation*} 
\Lambda_\beta(u)\left(\lambda(h_\beta)-u\frac{\Lambda_\beta'}{\Lambda_\beta} 
+u\frac{\tilde\Lambda_\beta'}{\tilde\Lambda_\beta}\right) =\tilde 
H_\beta\Lambda_\beta =0,
\end{equation*}
as operators on $V_\lambda^+$, so, since $\Lambda_\beta(u)$ is invertible,
\begin{equation*} \lambda(h_\beta)\Lambda_\beta(u)\tilde\Lambda_\beta(u) = 
u(\Lambda_\beta'\tilde\Lambda_\beta-\Lambda_\beta\tilde\Lambda_\beta').
\end{equation*}
Note that both sides of this equation make sense as power series in $u$ since 
$\Lambda_\beta(u)$ is already known by (i) to involve only finitely many 
positive powers of $u$. Hence, as series with only finitely many positive powers 
(but possibly infinitely many negative powers), we have
\begin{equation*}
\left(\frac{\Lambda_\beta}{\tilde\Lambda_\beta}\right)' 
=\lambda(h_\beta)u^{-1}\left(\frac{\Lambda_\beta}{\tilde{\Lambda}_\beta}\right),
\end{equation*}
and  so
\begin{equation*}
\frac{\Lambda_\beta(u)}{\tilde\Lambda_\beta(u)} = A_\beta 
u^{\lambda(h_\beta)},\end{equation*}
where $A_\beta$ is an operator on $V_\lambda^+$ independent of $u$. Equating 
coefficients of $u^{\lambda(h_\beta)}$ shows that $A_\beta 
=\Lambda_{\beta,\lambda(h_\beta)}$ and then the equation (of operators on 
$V_\lambda^+$) 
\begin{equation*} \Lambda_{\beta,\lambda(h_\beta)}(u)\tilde\Lambda_\beta(u) = 
u^{-\lambda(h_\beta)}\Lambda_\beta(u),\end{equation*}
proves both (iv) and (v).
\end{pf}

Fix a total order $\le$ on $R^+$.

\begin{prop}{\label{wint1}} 
Let $V$ be an integrable $\bu^e$-module, let $\lambda\in P^e_+$ and let $0\ne 
v\in V_\lambda^+$ 
be such that $V=\bu^e.v$.
\begin{enumerate}
\item[(i)]
If $V_\mu\ne 0$, then $\mu ={\lambda-\eta+r\delta}$ for some $\eta \in Q_+$, 
$r\in\bz$ such that $V^{fin}(\lambda)_{\lambda-\eta}\ne 0$. 

\item[(ii)] $V$ is spanned by the elements 
\begin{equation*} 
x_{\beta_1,r_1}^-x_{\beta_2, r_2}^-\cdots x_{\beta_s, r_s}^-
\bu(0).v,
\end{equation*}
for $0\le r_t<\lambda(h_{\beta_t})$, $0\le t\le s$, 
$\beta_1\le\beta_2\le\cdots\le\beta_s\in R^+$.
\end{enumerate} 
\end{prop}
\begin{pf} Since $\bu(>)_+.v =0$, it follows that $V =\bu(<)\bu(0).v$ and hence
\begin{equation*} 
V_\mu\ne 0\implies \mu=\lambda-\eta+r\delta,
\end{equation*}
for some $\eta\in Q^+$ and $r\in\bz$. Choose $\sigma\in W$ such that 
$\sigma(\lambda-\eta)\in P_+$. Since $V$ is integrable, $V_{\sigma(\mu)}\ne 0$, 
hence $\sigma(\lambda-\eta)=\lambda-\eta'$ for some $\eta'\in Q_+$. This implies 
that $V^{fin}(\lambda)_{\sigma(\lambda-\eta)}\ne 0$, and hence that 
$V^{fin}(\lambda)_{\lambda-\eta}\ne 0$.
  
To prove (ii), note first that it is clear that elements of the form
\begin{equation*}
x_{\beta_1,r_1}^-x_{\beta_2, r_2}^-\cdots x_{\beta_s, r_s}^- 
\bu(0).v,\ \ (0\le t\le s,\ \beta_t\in R^+, \ r_t\in \bz)
\end{equation*}
span $V$. We prove by induction on $s$ that any such element is in the span of 
the elements 
\begin{equation*} 
x_{\gamma_1,k_1}^-x_{\gamma_2, k_2}^-\cdots x_{\gamma_m, k_m}^- 
\bu(0).v,\  (0\le k_t<\lambda(h_{\gamma_t}),\ \ 
\gamma_1\le\gamma_2\le\cdots\le\gamma_m,\  0\le t\le m).
\end{equation*}

For $s=1$, $r_1\ge\lambda(h_{\beta_1})$, we have, by Proposition \ref{1.2}(iii), 
that
\begin{equation*}
\sum_{r=0}^{r_1}x^-_{\beta_1, r}\Lambda_{\beta_1, r_1-r}\bu(0).v 
=0.\end{equation*}
Since $\Lambda_{\beta_1,0}=1$, this implies that $x^-_{\beta_1, r_1}.v$ is in 
the span of the elements $x^-_{\beta_1, r}\bu(0).v$ with $0\le r<r_1$, from 
which the assertion follows. If $r_1<0$, we use
\begin{equation*}
\sum_{r=0}^{\lambda(h_{\beta_1})} x^-_{\beta_1, 
r+r_1}\Lambda_{\beta_1,\lambda(h_{\beta_1})-r}.v =0,
\end{equation*}
which follows from parts (i) and (ii) of Proposition 1.1. Since, by Proposition 
\ref{1.2}(v),  
$\Lambda_{\beta,\lambda(h_\beta)}$ is invertible, it follows that $x_{\beta_1, 
r_1}^-.v$ is in the span of the elements $x_{\beta_1, r}^-.v$ ($0\ge r> r_1$). 
An obvious induction now gives the result.

Suppose that we know the result for some $s\ge 1$. Let $\beta_t\in R^+$ and 
$r_t\in\bz$ for $0\le t\le s$ be such that $0\le r_t<\lambda(h_{\beta_t})$ for 
all $t$. We have
\begin{align*}
x^-_{\beta_0, r_0}x^-_{\beta_1,r_1}x^-_{\beta_2, r_2}\cdots 
x^-_{\beta_s, r_s}.v&=
x^-_{\beta_1, r_1}x^-_{\beta_0,r_0}x^-_{\beta_2, r_2}\cdots 
x^-_{\beta_s,r_s}.v\\
&\ \ \ \ \ +[x^-_{\beta_0, r_0}, x^-_{\beta_1,r_1}]x^-_{\beta_2, r_2}\cdots 
x^-_{\beta_s, r_s}.v.
\end{align*}
Since $[x^-_{\beta_0, r_0}, x^-_{\beta_1,r_1}]\in\frak g\otimes t^{r_0+r_1}$, 
the induction 
hypothesis applies to the second term on the right-hand side of 
the equality. As for the first term, notice that the induction hypothesis 
applied to $x^-_{\beta_0,r_0}x^-_{\beta_2, r_2}\cdots x^-_{\beta_s, r_s}.v$ 
implies that this is in the span of the elements obtained by applying ordered 
monomials 
in the $x^-_{\gamma, k}$ to $v$, for $\gamma\in R^+$ and $0\le 
k<\lambda(h_\gamma)$. 
Thus, we must show that every element of the form 
\begin{equation*}x_{\beta_1,r_1}^-x_{\gamma_1,k_1}^-x_{\gamma_2, k_2}^-\cdots 
x_{\gamma_m, k_m}^- 
\bu(0).v,\end{equation*}
where $\gamma_1\le\gamma_2\le \cdots\le \gamma_m$ and $0\le 
k_t<\lambda(h_{\gamma_t})$ for $1\le t\le m$, 
can be rewritten in the desired form. If $\beta_1\le\gamma_1$, there is nothing 
to prove. Otherwise, we have 
\begin{equation*} [x_{\beta_1, r_1}^-, x_{\gamma_1, k_1}^-] 
\in\frak{g}_{\beta_1+\gamma_1}\otimes t^{r_1+k_1}\end{equation*}
and $\lambda(h_{\beta_1+\gamma_1})\le r_1+k_1$.  
Using the induction hypothesis gives the result.
\end{pf}

{\section {Maximal Integrable and Maximal Finite-Dimensional Modules}} 

In this section we define, for every $\lambda\in P_+^e$  an  integrable 
$\ebu$-module $W(\lambda)$. Further, for any $n$-tuple ${\mbox{\boldmath$\pi$}} 
= (\pi_1(u),\pi_2(u),\cdots, \pi_n(u))$ of polynomials $\pi_i(u)$ in an 
indeterminate $u$ with 
constant term 1 and degree $\lambda(h_i)$, we define a finite-dimensional 
quotient $\bu$-module $W({\mbox{\boldmath$\pi$}})$ of $W(\lambda)$.

For $\lambda\in P_+^e$, let $I_\lambda$ be the left ideal in the subalgebra 
$\bu(<)\bu(0)\bu(\frak h)$ of $\ebu$ generated by the following elements:
\begin{align*} 
& h-\lambda(h)\ \  (h\in\frak h),\ \ \   \ \ \  \   \Lambda_{i,m} \ \  (i\in I,\ 
|m|>\lambda(h_i)),\\
&\Lambda_{i,\lambda(h_i)}\Lambda_{i,-m} -\Lambda_{i,\lambda(h_i) -m} \ \ 
(i\in I,\ 1\le m\le \lambda(h_i)),\\
&\left(\tilde X_i^-(u)\Lambda_i^+(u)\right)_m\bu(0)\ \ (i\in I,\  
m\in\bz),\\
&\left(X_{i,0}^-(u)^r\Lambda_i^+(u)\right)_m\bu(0)\ \ (i\in I,\ \  
r\ge 1,\ |m|>\lambda(h_i)).
\end{align*}

Let $\tilde I_\lambda$  be the left ideal in $\bu$ generated by $I_\lambda$   
and the $x_{i,m}^+$ for all $i\in I$, $m\in\bz$, and let $\tilde I^e_\lambda$ be 
the left ideal in $\ebu$ generated by $\tilde I_\lambda$ and $d-\lambda(d)$. Set
\begin{equation*} 
W(\lambda) =\ebu/\tilde I^e_\lambda = \bu/\tilde{I}_\lambda.
\end{equation*} 
Clearly, $W(\lambda)$ is a left $\bu^e$-module (and  a left $\bu$--module) 
through left multiplication. Let 
$w_\lambda$ be the image of $1$ in $W(\lambda)$. Then,
\begin{equation*} \bu(>)_+.w_\lambda =0,\ \ 
W(\lambda)=\bu^e.w_\lambda=\bu.w_\lambda.\end{equation*}

Since $\tilde{I}_\lambda\bu(0)\subset \tilde{I}_\lambda$, we can and do  regard 
$W(\lambda)$ as a right $\bu(0)$--module as well. For $\eta\in Q^+$, we set
\begin{equation*} W(\lambda)[\eta] =\bigoplus_{r\in\bz} 
W(\lambda)_{\lambda-\eta+r\delta}.\end{equation*}
Clearly, $W(\lambda)[\eta]$ is a right $\bu(0)$--module for all $\eta\in Q^+$ 
and we have 
\begin{equation*} W(\lambda)=\bigoplus_{\eta\in Q^+} 
W(\lambda)[\eta]\end{equation*}
as right $\bu(0)$--modules.

Let $I_\lambda(0) =I_\lambda\cap \bu(0)$. By the PBW theorem, it is easy to see 
that
\begin{equation*}\label{rightsym}
\bu(0)/I_{\lambda}(0)\cong\bc[\Lambda_{i,m}, \Lambda_{i,\lambda(h_i)}^{-1}: i\in 
I, 1\le m\le\lambda(h_i)].\end{equation*}
In particular,
\begin{equation*} W(\lambda)[0] \cong \bc[\Lambda_{i,m}, 
\Lambda_{i,\lambda(h_i)}^{-1}: i\in I, 1\le m\le\lambda(h_i)]\end{equation*}
as right $\bu(0)$--modules.
It follows immediately from Proposition \ref{wint1}(ii) that, for all $\eta\in 
Q^+$, $W(\lambda)[\eta]$ is a finitely-generated $\bu(0)$--module.

Next, let ${\mbox{\boldmath$\pi$}} =(\pi_1,\pi_2,\cdots 
,\pi_n)$ be an $n$-tuple of polynomials in an indeterminate $u$ with constant 
term $1$, and define,  for $a\in\bc$, an element 
$\lambda_{\mbox{\boldmath$\pi$},a}\in(\eh)^*$ by 
setting $\lambda_{\mbox{\boldmath$\pi$},a}(h_i)={\text{deg}}\,\pi_i$ ($i\in I$) 
and 
$\lambda_{\mbox{\boldmath$\pi$},a}(d)=a$. Set 
\begin{equation}{\label{pipm}}
\pi_i^+ =\pi_i, \ \ \pi_i^-(u)=
\frac{u^{{\text {deg}}\pi_i}\pi_i(u^{-1})}{u^{{\text 
{deg}}\pi_i}\pi_i(u^{-1})|_{u=0}}.
\end{equation} 
Let $I_\bpi(0)$ be the maximal ideal in $\bu(0)$ generated by 
\begin{equation*}
(\Lambda_i^\pm(u)-\pi_i^\pm (u))_s \ \ (i\in I,\ s\ge 0),
\end{equation*} 
and let $\bc_\bpi =\bu(0)/I_\bpi(0)$ be the one--dimensional $\bu(0)$-module.
 Set 
\begin{equation*}
W(\bpi)=W(\lambda_{\bpi,a})\otimes_{\bu(0)}\bc_\bpi.
\end{equation*}
 Then, $W(\bpi)$ is a left $\bu$-module (but not a $\bu^e$-module) with 
$x\in\bu$ acting as $x\otimes 1$. Let $w_{\mbox{\boldmath$\pi$}}$ be the image 
of $1$ in $W({\mbox{\boldmath$\pi$}})$. Note that 
$\Lambda_i^\pm(u).w_{\bpi}=\pi_i^\pm(u)w_{\bpi}$ ($i\in I$). The assignment 
$w_{\lambda_{\bpi,a}}\mapsto 
w_{\mbox{\boldmath$\pi$}}$ extends to a surjective $\bu$-module homomorphism  
$W(\lambda_{\bpi,a})\to W({\mbox{\boldmath$\pi$}})$.

Recall that the affine Lie algebra $\hat{\frak g}$ is an extension of $L^e(\frak 
g)$ by a 1-dimensional central subalgebra $\bc c$. Any representation of 
$L^e(\frak g)$ is a representation of $\hat{\frak g}$ by making $c$ act as zero. 
Set $h_0=[e_0^+,e_0^-]$, the bracket being evaluated in $\hat{\frak g}$. Then, 
$h_0=c-h_\theta$. 

The following result is proved in \cite{CP2}. 
\begin{lem}{\label{univint}} Let $V$ be a $\widehat{\frak g}$-module generated 
by an element 
$v\in V_\lambda$ ($\lambda\in P_+^e$) such that
\begin{align*} 
(e_i^+)^{-\lambda(h_i)+1}. v = 0\ \ &\ \ {\text{if}} \ \  
\lambda(h_i)\le 0,\\
(e_i^-)^{\lambda(h_i)+1}. v =0\ \ & \ \ {\text{if}} \ \ \lambda(h_i)\ge 0,
\end{align*}
for $i\in I$. Then, $V$ is an integrable $\widehat{\frak 
g}$-module.\hfill\qedsymbol
\end{lem}
  
\begin{thm}{\label{wint}} 
\begin{enumerate}
\item[(i)] Let $\lambda\in P^e_+$. Then, $W(\lambda)$ is an integrable 
$\ebu$-module. 
\item[(ii)] Let ${\mbox{\boldmath$\pi$}}$ be an $n$-tuple of polynomials with 
constant term one. Then, $W({\mbox{\boldmath$\pi$}})$ is a finite-dimensional 
$\bu$-module. 
\end{enumerate}
\end{thm}
\begin{pf} 

To prove (i), note that by Lemma \ref{univint} it suffices to  show that 
\begin{align*}  e_i^+. w_\lambda =0,\ \ & \ \  
(e_i^-)^{\lambda(h_i)+1}.w_\lambda =0\ \ (i\in I),\\
e_0^-.w_\lambda =0,\ \ &\ \ (e_0^+)^{\lambda(h_\theta)+1}.w_\lambda 
=0.\end{align*}
Suppose that $i\in I$. Then, $e_i^\pm =x_{i,0}^\pm$ and 
it follows from the definition of $W(\lambda)$ that 
$e_i^+.w_\lambda =0$. By Lemma \ref{gar}(ii),
\begin{equation*} (x_{i,0}^-)^{\lambda(h_i)+1}.w_\lambda 
=\left(X_{i,0}^-(u)^{\lambda(h_i)+1}\Lambda_i^+(u)\right)_{\lambda(h_i)+1}.
w_\lambda =0.\end{equation*} 
In particular, this proves that $\bu^{fin}.w_\lambda$ is a finite-dimensional 
$\frak g$-module.

Turning to the case $i=0$,  notice that $e_0^- =x_{\theta, -1}^+$  is a linear 
combination of products of the $x_{i,m}^+$ ($i\in I$, $m\in\bz$).
Hence, $x_{\theta,-1}^+.w_\lambda=0$. 

For any $m\ge 0$, let $w_m =(e_0^+)^{m+1}.w_\lambda$. Suppose that $w_m\ne 0$. 
Since $[e_0^+, x_{\beta,0}^-]=0$, it follows that \begin{equation*}\ug(<).w_m 
=e_0^{m+1}\ug(<).w_\lambda\end{equation*}  
is finite-dimensional. Since $W(\lambda)[\eta] =0$ for all but finitely many 
$\eta\in Q^+$, it follows that 
 \begin{equation*} W_m=\ug.w_m = \bu^{fin}(>)\bu(\frak 
h)\bu^{fin}(<).w_m\end{equation*} 
is a finite-dimensional $\frak g$-module.  
Hence, for all $\sigma\in W$ (the Weyl group of $\frak g$), we have 
\begin{equation*} 
(W_m)_{\sigma(\lambda-(m+1)\theta +(m+1)\delta) }\ne 0.
\end{equation*}
Choosing $\sigma$ so that $\sigma(\theta) = -\alpha_i$ for some $i\in I$, 
we get 
\begin{equation*}
W(\lambda)_{\sigma(\lambda)+(m+1)\alpha_i+(m+1)\delta}\ne 0. 
\end{equation*}
But this can only happen for finitely many values of $m$. 

This proves that $w_m =0$ for all but finitely many $m$.  The Lie subalgebra of 
$L(\frak g)$ generated by $e_0^\pm$ and $h_\theta$ is isomorphic to $sl_2$,  
and we have just shown that the corresponding $sl_2$-submodule generated by 
$w_\lambda$ is finite-dimensional. It follows from standard results that 
$(e_0^-)^{\lambda(h_\theta)+1}.w_\lambda =0$.

To prove (ii), it suffices now  to notice (using Proposition \ref{wint1}(ii)) 
that $W(\bpi)$ is spanned by the elements  
\begin{equation*} x_{\beta_1,r_1}^-x_{\beta_2, r_2}^-\cdots 
x_{\beta_s, r_s}^-\otimes 1\end{equation*}
for $s\ge 0$,  $0\le r_t<\lambda(h_{\beta_t})$, $\beta_t\in R^+$. 
 \end{pf}

The modules $W(\lambda)$ and $W({\mbox{\boldmath$\pi$}})$ have certain universal 
properties. 
\begin{prop}{\label{universal}} 

\vskip 12pt

\begin{enumerate}
\item[(i)] Let $V$ be any integrable  $\ebu$-module generated by a non-zero 
element $v\in V_\lambda^+$. 
 Then, $V$ is a quotient of $W(\lambda)$.
\item[(ii)] Let $V$ be a finite-dimensional quotient $\bu$-module  
of $W(\lambda)$, and assume that ${\text{dim}}\ V_\lambda =1$. Then, $V$ is a 
quotient of 
$W({\mbox{\boldmath$\pi$}})$ for some choice of ${\mbox{\boldmath$\pi$}}$.
\item[(iii)] Let $V$ be  finite-dimensional $\bu$-module generated by a vector 
$v\in V_\lambda^+$ and such that ${\text{dim}}\ V_\lambda =1$.   Then, $V$ is a 
quotient of $W({\mbox{\boldmath$\pi$}})$ for some  
${\mbox{\boldmath$\pi$}}$.\end{enumerate}
\end{prop}  

\begin{pf}
Part (i) is immediate from Proposition \ref{1.2} and the definition of 
$W(\lambda)$. 

To prove (ii), let $v\ne 0$ be the image of $w_\lambda$ in $V$. Notice that 
${\text{dim}}V_\lambda =1$ implies that
\begin{equation*}\Lambda_{\beta, m} .v =\pi_{\beta, m}v,\end{equation*}
for some scalars $\pi_{\beta, m}\in\bc$.  By Proposition \ref{1.2}(i), it 
follows that 
\begin{equation*} 
\pi_{\beta, m} = 0
\ \ \ \ \ (|m|>\lambda(h_\beta)).
\end{equation*}
For $i\in I$, set $\pi_i(u)=\sum_{k=0}^{\lambda(h_i)} \pi_{\alpha_i,k}u^k$. The 
$\pi_i(u)$ are 
polynomials 
with constant term $1$ and Proposition \ref{1.2}(v) shows that 
\begin{equation*}
\Lambda_i^\pm (u).v =\pi_i^\pm(u)v.
\end{equation*}
This shows that $V$ is a quotient of $W({\mbox{\boldmath$\pi$}})$, where 
${\mbox{\boldmath$\pi$}}$ is the $n$-tuple of 
polynomials defined above.

The proof of (iii) is identical.
\end{pf}

\vskip 24pt

\section{A tensor product theorem for $W({\mbox{\boldmath$\pi$}})$}
The main result of this section is the following theorem.

\begin{thm}{\label{wtensor}} 
Let ${\mbox{\boldmath$\pi$}}=(\pi_1(u),\pi_2(u),\cdots ,\pi_n(u))$ and 
${\mbox{\boldmath$\tilde\pi$}} =(\tilde\pi_1(u),\tilde\pi_2(u),\cdots 
,\tilde\pi_n(u))$ 
be $n$-tuples of polynomials in $u$ with constant term $1$, such that $\pi_i$ 
and $\tilde\pi_j$ are coprime for all $1\le i,j\le n$. Then,
\begin{equation*} 
W({\mbox{\boldmath$\pi$}}\ {\mbox{\boldmath$\tilde\pi$}}) \cong 
W({\mbox{\boldmath$\pi$}})\otimes W({\mbox{\boldmath$\tilde\pi$}})
\end{equation*} 
as $\bu$-modules, where ${\mbox{\boldmath$\pi$}}\ {\mbox{\boldmath$\tilde\pi$}}
=(\pi_1\tilde\pi_1,\pi_2\tilde\pi_2,\cdots ,\pi_n\tilde\pi_n)$.
\end{thm}

For $\beta\in R^+$ and $\bpi$ as in Theorem 2, define $\pi_\beta(u)$ by
\begin{equation}\label{pibeta}
\Lambda^+_\beta(u).w_{\mbox{\boldmath$\pi$}} = 
\pi_\beta(u)w_{\mbox{\boldmath$\pi$}}. 
\end{equation}

\begin{lem}{\label {pibeta1}} Let $\beta =\sum_{i=1}^n r_i\alpha_i\in R^+$. 
Then,
\begin{equation*} 
\pi_\beta(u) =\prod_{i=1}^n\pi_i^{r_id_i/d_\beta}(u).
\end{equation*} 
If $\theta_s\in R^+$ is the highest short root of $\frak g$, then $\pi_\beta$ 
divides $\pi_{\theta_s}$ for all $\beta\in R^+$.
\end{lem}
\begin{pf} The first statement is immediate from the formula for $h_\beta$ in 
terms of the $h_i$ given in Section 1 and the definition of the 
$\Lambda^+_\beta$. The second statement is proved case by case, using the 
explicit formula for $\theta_s$ (there is a simple uniform proof when $\frak g$ 
is simply-laced).
\end{pf}
 
From now on, we shall assume that, given $\beta\in R^+$ and an $n$-tuple of 
polynomials ${\mbox{\boldmath$\pi$}}$, we have defined a polynomial 
$\pi_\beta(u)$ by the formula given in (\ref{pibeta}). The polynomials 
$\pi_\beta^\pm(u)$ are defined as in (\ref{pipm}).

\begin{defn} Let ${\mbox{\boldmath$\pi$}}=\{\pi_1,\cdots ,\pi_n\}$ be an 
$n$-tuple of polynomials in $u$ with constant term 1. Define 
$M({\mbox{\boldmath$\pi$}})$ to be the left $\bu$-module obtained by taking the 
quotient of $\bu$ by the left ideal generated by the following:
\begin{align*}
& x^+_{i,k},\ \  \ \  \ h_{i,0}-{\text{deg}}\pi_i\ \ (i\in I, k\in\bz),\\
&\left(\Lambda_i^\pm(u) -\pi_i^\pm(u)\right)_s\ \ (i\in I, s\ge 0),\\
& \left(\pi_\beta(u)\tilde X_\beta^-(u)\right)_s\ \ (\beta\in R^+, s\in\bz).
\end{align*}
\end{defn}
Let $m_{\mbox{\boldmath$\pi$}}$ be the image of 1 in 
$M({\mbox{\boldmath$\pi$}})$.  It is clear from equation (\ref{pibeta}) and 
Proposition 1.1(iii) that, for all $\beta\in R^+$ and all $s\in\bz$, 
\begin{align*}
\left(\Lambda_\beta^\pm(u)-\pi_\beta^\pm(u)\right)_s.m_{\mbox{\boldmath$\pi$}}
&=0,\\
\left(\pi_\beta(u)\tilde H_\beta(u)\right)_s.m_\bpi& =0.
\end{align*}
Set $\lambda_\bpi(h_i)=\text{deg}\,\pi_i$.  
\begin{lem}{\label{mdimw}} We have
\begin{equation*} M(\bpi) =\bigoplus_{\eta\in Q^+} 
M(\bpi)_{\lambda_\bpi-\eta},\end{equation*}
and ${\text{dim}}\, M(\bpi)_{\lambda-\eta}<\infty$. Further,
for all $\beta\in R^+$, $s\in\bz$, we have
\begin{equation*}\left(\pi_{\theta_s}(u)\tilde{X}_\beta^-(u)\right)_s.M(\bpi) 
=0.\end{equation*}
\end{lem}
\begin{pf} The set
\begin{equation*} 
\{x_{\beta, r}^-:\beta\in R^+, 0\le 
r<\lambda(h_\beta)\}\cup\{x^-_\beta\otimes t^m\pi_\beta(t): \beta\in R^+, 
m\in\bz\}
\end{equation*}
is a basis of 
\begin{equation*} 
\left(\bigoplus_{\beta\in R^+} \bc x_\beta^-\right)\otimes 
\bc[t,t^{-1}].
\end{equation*}
By the PBW basis theorem, we can write
\begin{equation*}\bu(<) 
=\bu_{\mbox{\boldmath$\pi$}}\bu^{\mbox{\boldmath$\pi$}}\end{equation*}
where $\bu_{\mbox{\boldmath$\pi$}}$ (resp. $\bu^{\mbox{\boldmath$\pi$}}$)  
consists of ordered monomials from the first 
(resp. second) set. The relation
\begin{equation*} \left(\pi_\beta(u)\tilde X_\beta^-(u)\right)_s. 
m_{\mbox{\boldmath$\pi$}} 
=0\end{equation*}
for all $s\in\bz$ implies that 
$(\bu^{\mbox{\boldmath$\pi$}})_+.m_{\mbox{\boldmath$\pi$}} =0$. A further use of 
the PBW theorem now shows that 
\begin{equation*} M({\mbox{\boldmath$\pi$}})\cong \bu_{\mbox{\boldmath$\pi$}} 
\end{equation*}
as vector spaces. Moreover, this isomorphism takes 
$M(\bpi)_{\lambda_\bpi-\eta}$ to
\begin{equation*} \bu_\bpi(\eta) =\{x\in\bu_\bpi: [h,x]=\eta(h)x\ \ \forall 
h\in\frak h\}.\end{equation*}
Since this space is clearly finite--dimensional, the first statement of the 
lemma follows. 

For the second statement, we show by induction on ${\text{ht}}\,\eta$ that
\begin{equation*}\left(\pi_\theta(u)\tilde{X}_\beta^-(u)\right)_s.
M(\bpi)_{\lambda_\bpi-\eta} =0\ \ \ \ \ (\beta\in R^+,\ s\in\bz).
\end{equation*}
If $\eta =0$, the result is immediate from the definition of $M(\bpi)$ and the 
last part of Lemma 3.1. In general, let $x_{\beta_1, r_1}^-\cdots x_{\beta_k, 
r_k}^-.m_\bpi\in M(\bpi)_{\lambda_\bpi-\eta}$.
Then,
\begin{align*}  
\left(\pi_{\theta_s}(u)\tilde{X}_\beta^-(u)\right)_sx_{\beta_1, r_1}^-\cdots 
x_{\beta_k, r_k}^-.m_\bpi&= x_{\beta_1, r_1}^- 
\left(\pi_\theta(u)\tilde{X}_\beta^-(u)\right)_sx_{\beta_2,r_2}^-\cdots 
x_{\beta_k, r_k}^-.m_\bpi\\
&+\left(\pi_\theta(u)\tilde{X}_{\beta+\beta_1}^-(u)\right)_{s+r_1}x_{\beta_2,r_2
}^-\cdots x_{\beta_k, r_k}^-.m_\bpi,\end{align*}
where we understand that $\tilde{X}_{\beta+\beta_1}^-(u) =0$ if 
$\beta+\beta_1\notin R^+$.
The right-hand side is zero by the induction hypothesis, and the inductive step 
is complete.
\end{pf}
\begin{lem}
The $\bu$-module $W({\mbox{\boldmath$\pi$}})$ is a quotient of 
$M({\mbox{\boldmath$\pi$}})$, and any finite-dimensional 
quotient of $M({\mbox{\boldmath$\pi$}})$ is a quotient of 
$W({\mbox{\boldmath$\pi$}})$.\end{lem}
\begin{pf}Let $V$ be a finite-dimensional quotient of $M(\bpi)$, let $v\in V$ be 
the image of $m_\bpi$, and let $\lambda=\lambda_\bpi$. Then, 
${\text{dim}}\,V_\lambda={\text{dim}}\,M(\bpi)_\lambda=1$, so by Proposition 
2.1(ii), $V$ is a quotient of some $W(\tilde\bpi)$ with 
$\lambda=\lambda_{\tilde\bpi}$. Since ${\text{dim}}\,W(\tilde\bpi)_\lambda=1$, 
$v$ is a scalar multiple of the image of $w_{\tilde\bpi}$. But then by comparing 
the action of $\Lambda_i^\pm(u)$ on $w_{\tilde\bpi}$ and on $m_\bpi$, we see 
that $\bpi=\tilde\bpi$.

To show that $W({\mbox{\boldmath$\pi$}})$ is a 
quotient of $M({\mbox{\boldmath$\pi$}})$, it is clear from the definitions of 
these modules that we only need to show that 
\begin{equation*}\left(\pi_\beta(u)\tilde 
X_\beta^-(u)\right)_s.w_{\mbox{\boldmath$\pi$}} =0\ \ \ \ (s\in\bz).
\end{equation*}
Since  $W({\mbox{\boldmath$\pi$}})$ is a quotient of 
$W(\lambda_{\mbox{\boldmath$\pi$},0})$, this follows from Proposition \ref{1.2},
thus completing the proof of the lemma.
\end{pf}

Denote by $\Delta:\bu\to\bu\otimes\bu$ the comultiplication of $\bu$ defined by 
extending to an algebra homomorphism the assignment 
\begin{equation*} 
x\to x\otimes 1+ 1\otimes x,
\end{equation*}
for all $x\in \loopg$. The following lemma is proved in \cite{G}.
\begin{lem}{ \label{comult}} For all $\beta\in R^+$,
\begin{equation*} 
\Delta(\Lambda_\beta^\pm) =\Lambda^\pm_\beta\otimes\Lambda^\pm_\beta,
\end{equation*}
where
\begin{equation*} 
\Lambda^\pm_\beta\otimes \Lambda^\pm_\beta  =\sum_{k,m\ge 0} 
(\Lambda_{\beta, \pm k}\otimes\Lambda_{\beta, \pm m})u^{k+m}.
\end{equation*}
\hfill\qedsymbol\end{lem}

Theorem \ref{wtensor} is now clearly a consequence of the following  
proposition.

\begin{prop}{\label{mtensor}} Assume that 
${\mbox{\boldmath$\pi$}}=\{\pi_1,\pi_2,\cdots ,\pi_n\}$ 
and ${\mbox{\boldmath$\tilde\pi$}} =\{\tilde\pi_1,\tilde\pi_2,\cdots 
,\tilde\pi_n\}$ are $n$-tuples 
of polynomials with constant term $1$, such that $\pi_i$ and $\tilde\pi_j$ are 
coprime for all $1\le i,j\le n$. Then:
\begin{enumerate}
\item[(i)] $M({\mbox{\boldmath$\pi$}}\ {\mbox{\boldmath$\tilde \pi$}}) \cong 
M({\mbox{\boldmath$\pi$}})\otimes M({\mbox{\boldmath$\tilde\pi$}})$;
\item[(ii)] every finite-dimensional quotient $\bu$-module of 
$M({\mbox{\boldmath$\pi$}})\otimes M({\mbox{\boldmath$\tilde\pi$}})$ is a 
quotient of $W(\bpi)\otimes W(\tilde\bpi)$. 
\end{enumerate}
\end{prop}
\begin{pf} Set $\lambda =\lambda_\bpi+\lambda_{\tilde\bpi}$.
 It is clear from the proof of Lemma 
\ref{mdimw} that, for all $\eta\in Q_+$, we have 
\begin{equation*} 
M({\mbox{\boldmath$\pi$}}\ {\mbox{\boldmath$\tilde\pi$}})_{\lambda-\eta} \cong 
\left(M({\mbox{\boldmath$\pi$}})\otimes 
M({\mbox{\boldmath$\tilde\pi$}})\right)_{\lambda-\eta}
\end{equation*}
{\it as (finite-dimensional) vector spaces}. 
To prove (i), it therefore suffices to prove that there 
exists a surjective homomorphism of $\bu$-modules $M({\mbox{\boldmath$\pi$}}\ 
{\mbox{\boldmath$\tilde\pi$}}) \to 
M({\mbox{\boldmath$\pi$}})\otimes M({\mbox{\boldmath$\tilde\pi$}})$. It is easy 
to 
see, using Lemma \ref{comult}, that the 
element $m_{\mbox{\boldmath$\pi$}}\otimes m_{\mbox{\boldmath$\tilde\pi$}}$ 
satisfies 
the defining relations of
$M({{\mbox{\boldmath$\pi$}}{\mbox{\boldmath$\tilde\pi$}}})$, so there exists  a 
$\bu$-module map $M({\mbox{\boldmath$\pi$}}\ {\mbox{\boldmath$\tilde\pi$}}) \to 
M({\mbox{\boldmath$\pi$}})\otimes M({\mbox{\boldmath$\tilde\pi$}})$ that sends 
$m_{{\mbox{\boldmath$\pi$}}{\mbox{\boldmath$\tilde\pi$}}}$ to 
$m_{\mbox{\boldmath$\pi$}}\otimes m_{\mbox{\boldmath$\tilde\pi$}}$.  Thus, 
to prove (i), we must show that, if $\pi_i$ and $\tilde\pi_j$ have no roots in 
common, the element $m_{\mbox{\boldmath$\pi$}}\otimes 
m_{\mbox{\boldmath$\tilde\pi$}}$  
generates $M({\mbox{\boldmath$\pi$}})\otimes M({\mbox{\boldmath$\tilde\pi$}})$ 
as a $\bu$-module.

Set \begin{equation*} N =\bu.(m_{\mbox{\boldmath$\pi$}}\otimes 
m_{{\mbox{\boldmath$\tilde\pi$}}}).\end{equation*} 
Assume that, for 
all $\eta =\sum_i r_i\alpha_i$, $\tilde\eta=\sum_i \tilde r_i\alpha_i$, with 
${\text{ht}}\,\eta = \sum_ir_i< s$, ${\text{ht}}\,\tilde\eta =\sum_i\tilde 
r_i<s$, we have 
\begin{equation*}M({\mbox{\boldmath$\pi$}})_{\lambda_{\mbox{\boldmath$\pi$}}
-\eta}\otimes 
M({\mbox{\boldmath$\tilde 
\pi$}})_{\lambda_{{\mbox{\boldmath$\tilde\pi$}}}-\tilde\eta}\subset 
N.\end{equation*}
We shall prove that 
\begin{align*} 
&(x_{i,m}^-.M({\mbox{\boldmath$\pi$}})_{\lambda_{\mbox{\boldmath$\pi$}}-\eta})
\otimes M({\mbox{\boldmath$\tilde 
\pi$}})_{\lambda_{{\mbox{\boldmath$\tilde\pi$}}}-\tilde\eta}
\subset N,\\
&M({\mbox{\boldmath$\pi$}})_{\lambda_{{\mbox{\boldmath$\pi$}}}-\eta}\otimes 
x_{i,m}^-.(M({\mbox{\boldmath$\tilde\pi$}})_{\lambda_{{\mbox{\boldmath$\tilde\pi
$}}}-
\tilde \eta})\subset N
\end{align*}
for all $i\in I$, $m\in\bz$. This will prove that 
\begin{equation*}
M({\mbox{\boldmath$\pi$}})_{\lambda_{\mbox{\boldmath$\pi$}}-\eta}\otimes 
M({\mbox{\boldmath$\tilde\pi$}})_{\lambda_{{\mbox{\boldmath$\tilde\pi$}}}-\tilde
\eta}\subset N,\end{equation*}
when ${\text{ht}}\,\eta\le s$, 
${\text{ht}}\,\tilde\eta\le s$, and hence, by induction on $s$, that 
$N=M({\mbox{\boldmath$\pi$}})\otimes 
M({\mbox{\boldmath$\tilde\pi$}})$.

Since $\pi_{\theta_s}$ and $\tilde\pi_{\theta_s}$ are coprime, we can choose 
polynomials $R(u)$, $\tilde R(u)$ such that
\begin{equation*} 
R\pi_{\theta_s}+\tilde R\tilde\pi_{\theta_s} =1.
\end{equation*}
By the second part of Lemma 3.2,
\begin{align*} 
\left(R\pi_{\theta_s}\tilde X_i^-(u)\right)_m.w &=0,\\ 
\left(\tilde R\tilde\pi_{\theta_s}\tilde X_i^-(u)\right)_m.\tilde w &=0,
\end{align*}
for all $i\in I$, $m\in\bz$, $w\in M(\bpi)$, $\tilde w\in M(\tilde\bpi)$. 
Hence,
\begin{align*}
\left(R\pi_{\theta_s}\tilde X_i^-(u)\right)_m & .(w\otimes \tilde w)\\ 
&= \left(R\pi_{\theta_s}\tilde X_i^-(u)\right)_m. w\otimes \tilde w + w\otimes 
\left(R\pi_{\theta_s}\tilde X_i^-(u)\right)_m. \tilde w\\
&=  w\otimes \left((1-\tilde R\tilde\pi_{\theta_s})\tilde 
X_i^-(u)\right)_m.\tilde w\\
&= w\otimes x_{i,m}^-.\tilde w.\end{align*}
Taking $w\in M({\mbox{\boldmath$\pi$}})_{\lambda_{\mbox{\boldmath$\pi$}}-\eta}$ 
and $\tilde w\in M({\mbox{\boldmath$\tilde 
\pi$}})_{\lambda_{\mbox{\boldmath$\tilde \pi$}}-\eta}$, 
so that $w\otimes \tilde w\in N$,   it follows that 
$w\otimes x_{i,m}^-.\tilde w\in N$ for all $i\in I$, $m\in\bz$. The other 
inclusion is proved similarly, and the 
proof of part (i) is complete.

Suppose that $V$ is a finite-dimensional quotient of 
$M({\mbox{\boldmath$\pi$}})\otimes M({\mbox{\boldmath$\tilde\pi$}})$ 
with kernel $K$. We shall prove that, for all $r\ge 1$, $i\in I$, 
$s\ge\lambda_{\mbox{\boldmath$\pi$}}(h_i)+1$, $\tilde s\ge 
\lambda_{{\mbox{\boldmath$\tilde\pi$}}}(h_i)+1$,
\begin{align}\label{*}
\left(\pi_i(u)X_i^-(u)^r\right)_s.m_{\mbox{\boldmath$\pi$}}\otimes 
M({{\mbox{\boldmath$\tilde\pi$}}})\subset K,&\\ 
M({\mbox{\boldmath$\pi$}})\otimes 
\left(\pi_i(u)X_i^-(u)^r\right)_s.m_{{\mbox{\boldmath$\tilde\pi$}}}\subset K. &
\end{align}
Since the sum of these subspaces is the kernel of the quotient map
\begin{equation*}
M({\mbox{\boldmath$\pi$}})\otimes M({\mbox{\boldmath$\tilde\pi$}})\to 
W({\mbox{\boldmath$\pi$}})\otimes W({\mbox{\boldmath$\tilde\pi$}}),
\end{equation*}
it follows that $V$ is a quotient of $W({\mbox{\boldmath$\pi$}})\otimes 
W({\mbox{\boldmath$\tilde\pi$}})$, which proves part 
(ii).

To prove equation (\ref{*}), it suffices  (by the proof of Proposition 
\ref{1.2}(ii))  to prove that, for all $i\in I$, $m\in\bz$,
\begin{equation*} 
(x_{i,m}^-)^{r_i+1}.m_{\mbox{\boldmath$\pi$}}\otimes 
m_{{\mbox{\boldmath$\tilde\pi$}}}\in K,
\end{equation*}
where $r_i ={\text{deg}}\,\pi_i=\lambda_\bpi(h_i)$.
Since $V$ is finite-dimensional,  
the element $(x_{i,m}^-)^{r}.m_{\mbox{\boldmath$\pi$}}\otimes 
m_{{\mbox{\boldmath$\tilde\pi$}}}\in K$ for some $r\ge 0$. Let $r_0$ be the 
smallest value of $r$ with this property. Since
\begin{equation*} 
x_{i,-m}^+(x_{i,m}^-)^{r_0}.m_{\mbox{\boldmath$\pi$}}\otimes 
m_{{\mbox{\boldmath$\tilde\pi$}}} =(r_i-r_0+1) 
(x_{i,m}^-)^{r_0-1}.m_{\mbox{\boldmath$\pi$}}\otimes 
m_{{\mbox{\boldmath$\tilde\pi$}}},
\end{equation*}
it follows by the minimality of $r_0$ that $r_i +1 =r_0$. 

Equation (3.3) is proved similarly, and we are done.\end{pf}
\medskip 

Note that, since 
${\text{dim}}\,W(\bpi)_{\lambda_\bpi} =1$, it follows that $W(\bpi)$ has a 
unique irreducible quotient $V(\bpi)$. Write $\bpi$ as a product
\begin{equation*}
\bpi =\bpi^{(1)}\bpi^{(2)}\cdots \bpi^{(k)},
\end{equation*}
where $\bpi^{(j)}$ is such that
\begin{equation*}
\pi_{\theta_s}^{(j)}= (1-a_ju)^{m_j},
\end{equation*}
for some $m_j>0$ and $a_j\in\bc^\times$ with $a_j\ne a_k$ if $j\ne k$.
The following result was proved in \cite{CP1}. 
\begin{prop}{\label{vtensor}} With the above notation,
\begin{equation*} 
V(\bpi)\cong V(\bpi^{(1)})\otimes V(\bpi^{(2)})\otimes\cdots 
\otimes V(\bpi^{(k)})
\end{equation*} 
as $L(\frak g)$-modules. Further, 
\begin{equation*} 
V(\bpi^{(j)})\cong V^{fin}(\lambda_{\bpi_j})
\end{equation*}
as $\frak g$-modules.
\hfill\qedsymbol\end{prop}

\section{The quantum case} 
In the remainder of this paper, we shall assume that $\frak g$ is simply-laced.
 
Let $q$ be an indeterminate, let $\bc(q)$ be the field of rational
functions in $q$ with complex coefficients, and let $\ba=\bc[q,q^{-1}]$ be
the subring of Laurent polynomials. For $r,m\in\bn$, $m\ge r$, define
\begin{equation*} 
[m]=\frac{q^m -q^{-m}}{q -q^{-1}},\ \ \ \ [m]! =[m][m-1]\ldots [2][1],\ \ \ \ 
\left[\begin{matrix} m\\ r\end{matrix}\right] 
= \frac{[m]!}{[r]![m-r]!}.
\end{equation*}
Then, $\left[\begin{matrix} m\\r\end{matrix}\right]\in\ba$. 

Let $\bu^e_q$ be the  quantized enveloping algebra over $\bc(q)$  associated to 
$\eloopg$. Thus, $\bu^e_q$ is the quotient of the quantum affine algebra 
obtained by setting the central generator equal to $1$. It follows from 
\cite{Dr2}, \cite{B}, \cite{J} that $\bu^e_q$ is the algebra 
with generators $\bx_{i,r}^{{}\pm{}}$ ($i\in I$, $r\in\bz$), $K_i^{{}\pm 1}$ 
($i\in I$), $\bh_{i,r}$ ($i\in I$, $r\in \bz\backslash\{0\}$), $D^{\pm 1}$, 
and the following defining relations:
\begin{align*}
   K_iK_i^{-1} = K_i^{-1}K_i
  =1, \ \ &\ \ DD^{-1} =D^{-1}D =1, \\
 K_iK_j =K_jK_i,\;\;
  &K_i\bh_{j,r} =\bh_{j,r}K_i,\\ 
 K_i\bx_{j,r}^\pm K_i^{-1} = q^{{}\pm
    a_{ij}}\bx_{j,r}^{{}\pm{}},\ \ & D\bx_{j,r}^\pm D^{ -1} =q^r\bx^\pm_{j,r},\\ 
  [\bh_{i,r},\bh_{j,s}]=0,\; \; & [\bh_{i,r} , \bx_{j,s}^{{}\pm{}}] =
  \pm\frac1r[ra_{ij}]\bx_{j,r+s}^{{}\pm{}},\\ 
 \bx_{i,r+1}^{{}\pm{}}\bx_{j,s}^{{}\pm{}} -q^{{}\pm
    a_{ij}}\bx_{j,s}^{{}\pm{}}\bx_{i,r+1}^{{}\pm{}} &=q^{{}\pm
    a_{ij}}\bx_{i,r}^{{}\pm{}}\bx_{j,s+1}^{{}\pm{}}
  -\bx_{j,s+1}^{{}\pm{}}\bx_{i,r}^{{}\pm{}},\\ [\bx_{i,r}^+ ,
  \bx_{j,s}^-]=\delta_{i,j} & \frac{ \psi_{i,r+s}^+ -
    \psi_{i,r+s}^-}{q - q^{-1}},\\ 
\sum_{\pi\in\Sigma_m}\sum_{k=0}^m(-1)^k\left[\begin{matrix}m\\k\end{matrix}
\right]
  \bx_{i, r_{\pi(1)}}^{{}\pm{}}\ldots \bx_{i,r_{\pi(k)}}^{{}\pm{}} &
  \bx_{j,s}^{{}\pm{}} \bx_{i, r_{\pi(k+1)}}^{{}\pm{}}\ldots
  \bx_{i,r_{\pi(m)}}^{{}\pm{}} =0,\ \ \text{if $i\ne j$},
\end{align*}
for all sequences of integers $r_1,\ldots, r_m$, where $m =1-a_{ij}$, $\Sigma_m$ 
is the symmetric group on $m$ letters, and the $\psi_{i,r}^{{}\pm{}}$ are 
determined by equating powers of $u$ in the formal power series 
$$\sum_{r=0}^{\infty}\psi_{i,\pm r}^{{}\pm{}}u^{{}\pm r} = K_i^{{}\pm 1} 
exp\left(\pm(q-q^{-1})\sum_{s=1}^{\infty}\bh_{i,\pm s} u^{{}\pm s}\right).$$ 

\vskip 12pt

Define the $q$-divided powers
\begin{equation*}(\bx_{i,k}^\pm)^{(r)} 
=\frac{(\bx_{i,k}^\pm)^r}{[r]!},\end{equation*}
for all $i\in I$, $k\in\bz$, $r\ge 0$.

Suppose that $a_{ij} =-1$. Then, a special case  of the above 
relations is
\begin{equation*} 
(\bx_{i,s}^-)^2\bx_{j,r}^- 
-(q+q^{-1})\bx_{i,s}^-\bx_{j,r}^-\bx_{i,s}^-+\bx_{j,r}^-(\bx_{i,s}^-)^2 = 0.
\end{equation*}
Set $\gamma_{s,r}^{i,j} =\bx_{i,s}^-\bx_{j,r}^- -q\bx_{j,r}^-\bx_{i,s}^-$.
Again, the relations in $\bu_q^e$ imply that
\begin{equation*}
\gamma_{s,r}^{i,j} = -\gamma^{j,i}_{r+1,s-1}.
\end{equation*}
The following result is proved in \cite{L2}.

\begin{prop}\label{ij} For $i,j\in I$ with $a_{ij}=-1$, $r,s,l,m\in\bz$, and 
$l,m\ge 0$, there exist $f_p\in\ba$, for $0\le p\le{\text{min}}(l,m)$, such that 
\begin{equation*}
(\bx_{i,s}^-)^{(l)}(\bx_{j,r}^-)^{(m)}=\sum_p 
f_p(\bx_{j,r}^-)^{(m-p)}(\gamma^{i,j}_{s,r})^{(p)}(\bx_{i,s}^-)^{(l-p)}.
\end{equation*}
Further, there exists $g_p\in\ba$, for $0\le p\le m$, such that 
\begin{equation*} 
(\gamma_{s,r}^{i,j})^{(m)}=\sum_p g_p(\bx_{j,r}^-)^{(p)} 
(\bx_{i,s}^-)^{(m)}(\bx_{j,r}^-)^{(m-p)}.\ \ \ \ \ \ \ \ \ \end{equation*}
\hfill\qedsymbol
\end{prop}

Define
\begin{equation*}
\blambda_i^\pm(u) =\sum_{m=0}^\infty {\mbox{\boldmath 
$\Lambda$}}_{i,\pm m}u^m ={\text{exp}}\left(-\sum_{k=1}^\infty \frac{\bh_{i,\pm 
k}}{[k]}u^k\right).
\end{equation*}
The subalgebras $\bu_q$, $\bu_q^{fin}$, $\bu_q(<)$, etc., are defined in the 
obvious way. Let $\bu_q(\eh) $ be the subalgebra of $\bu_q^e$ generated by 
$K_i^{\pm 1}$ ($i\in I$) and $D^{\pm 1}$. Let $\bu_q(0)$ be the subalgebra of 
$\bu_q^e$ generated by the 
elements $\blambda_{i,m}$ ($i\in I$, 
$m\in\bz$). The following result is a simple corollary of the PBW theorem for 
$\bu^e_q$ \cite{B}.

\begin{lem}\label{triangle} $\bu^e_q\cong 
\bu_q(<)\bu_q(0)\bu_q(\eh)\bu_q(>)$.\hfill\qedsymbol 
\end{lem}

For any invertible element $x\in\bu_q^e$, define
\begin{equation*}\left[\begin{matrix}x\\r\end{matrix}\right] = \frac{xq^{r} 
-x^{-1}q^{-r}}{q-q^{-1}}.\end{equation*}

Let $\bu^e_\ba$ be the $\ba$-subalgebra of $\bu^e_q$ generated by the $K_i^{\pm 
1}$, $(\bx_{i,k}^\pm)^{(r)}$  ($i\in I$, $k\in\bz$,  $r\ge 0$),  
$D^{\pm 1}$, and  $\left[\begin{matrix}D\\ r\end{matrix}\right]$ ($r\ge 1$). 
Then, \cite{L2}, \cite{BCP},
\begin{equation*}
\bu^e_q=\bc(q)\otimes_\ba\bu^e_\ba .
\end{equation*}

 Define $\bu_\ba(<)$,  $\bu_\ba(0)$ and $\bu_\ba(>)$ in the obvious way. 
Let $\bu_\ba(\eh)$ be the $\ba$--subalgebra of $\bu_\ba$ generated by the 
elements $K_i^{\pm 1}$, $D^{\pm 1}$, $\left[\begin{matrix}K_i\\ 
r\end{matrix}\right]$ and  $\left[\begin{matrix}D\\ r\end{matrix}\right]$ ($i\in 
I$, $r\in\bz$).
 The following is proved as in Proposition 2.7 in \cite{BCP}.
\begin{prop}
$\bu^e_\ba =\bu_\ba(<)\bu_\ba(0)\bu_\ba^e(\frak h)\bu_\ba(>)$.
\hfill\qedsymbol
\end{prop}

The next lemma is easily checked.
\begin{lem}
\begin{enumerate}
\item[(i)]  There is a unique $\bc$--linear  anti--automorphism $\Psi$ of 
$\bu_q^e$ such that $\Psi(q) =q^{-1}$ and
\begin{align*}\Psi(K_i) =K_i,\ \ &\ \Psi(D) =D,\\ \Psi(x_{i,r}^\pm) 
=x_{i,r}^\pm, \ \ &\ \ \Psi(h_{i,r}) = -h_{i,r},\end{align*}
for all $i\in I$, $r\in\bz$.
\item[(ii)] There is a unique algebra automorphism $\Phi$ of $\bu^e_q$ over 
$\bc(q)$ such that $\Phi(\bx_{i,r}^\pm)=\bx_{i,-r}^\pm$, 
$\Phi(\Lambda_i^+(u))=\Lambda_i^-(u)$.
\end{enumerate}
\end{lem}

The first part of the following  result is proved  in \cite{BCP}, and the second 
part follows from it by applying $\Psi$.

\begin{lem}\label{order} 
\begin{enumerate} 
\item[(i)] Let $s>s'$, $m,m'\ge 0$, $i\in I$. Then, 
$(\bx_{i,s}^-)^{(m)}(\bx_{i,s'}^-)^{(m')}$ is in the span of the 
elements
\begin{equation*}
(\bx_{i,r_1}^-)^{(k_1)}(\bx_{i,r_2}^-)^{(k_2)}
\cdots (\bx_{i,r_l}^-)^{(k_l)},
\end{equation*}
where $s'\le r_1<r_2<\cdots < r_l\le s$, $\sum_p k_p=m+m'$ and $\sum_p
k_pr_p=ms+m's'$.
\item[(ii)] Let $s<s'$, $m,m'\ge 0$, $i\in I$. Then,  
$(\bx_{i,s}^-)^{(m)}(\bx_{i,s'}^-)^{(m')}$ is in the span of the 
elements
\begin{equation*}
(\bx_{i,r_1}^-)^{(k_1)}(\bx_{i,r_2}^-)^{(k_2)}
\cdots (\bx_{i,r_l}^-)^{(k_l)},
\end{equation*}
where $s'\ge r_1>r_2>\cdots >r_l\ge s$, $\sum k_p=m+m'$ and $\sum 
k_pr_p=ms+m's'$.
\end{enumerate}
\hfill\qedsymbol
\end{lem}

For $i\in I$, let $\tilde\bX_i^-(u)$, $\bX_i^-(u)$, $\bX_{i,0}^-(u)$ be the 
formal power series with coefficients in $\bu_q$ given  by the same formulas as  
the 
$\tilde X_i^-(u)$, etc., in Section 1. The next result is a reformulation of 
results in \cite[Section 5]{CP6}.
\begin{lem}{\label{gi}}  Let $s\ge r\ge 1$, $i\in I$. 
\begin{enumerate}
\item[(i)] 
\begin{equation*}(\bx_{i,0}^+)^{(r)} (\bx_{i,1}^-)^{(s)} = 
(-1)^r(\bX_i^-(u)^{(s-r)}\blambda_i^+(u))_s\mod \bu_q\bu_q(>)_+.
\end{equation*}
\item[(ii)] 
\begin{equation*}(\bx_{i,1}^+)^{(r)}(\bx_{i,0}^-)^{(s)} = 
(-1)^r(\bX_{i,0}^-(u)^{(s-r)}\blambda_i^+(u))_s\mod \bu_q\bu_q(>)_+.
\end{equation*}
\item[(iii)] \begin{equation*}
\left(\bX_i^-(u)^{(r)}\right)_{s+r}=\sum\mu(s_0,s_1,\cdots)(\bx_{i,0}^-)^{(s_0)}
(\bx_{i,1}^-)^{(s_1)}\cdots ,\end{equation*}
where the sum is over those non-negative integers $s_0,s_1,\cdots $ such that 
$\sum_ts_t=r$ and $\sum_t ts_t =s+r$, and the coefficients $\mu(s_0,s_1,\cdots 
)\in\ba$. In particular, the coefficient of $(\bx_{i,s+1}^-)^{(r)}$ in 
$\left(\bX^-_i(u)^{(r)}\right)_{(s+1)r}$  is $q^{sr(r-1)}$.\hfill\qedsymbol
\end{enumerate}
\end{lem}

\begin{defn} A $\bu_q^e$--module $V_q$  is said to be of type 1 if 
\begin{equation*}
V_q=\bigoplus_{\lambda\in P^e}(V_q)_\lambda,
\end{equation*}
where
\begin{equation*} 
(V_q)_\lambda =\{v\in V_q: K_i.v =q^{\lambda(h_i)}v\ \forall 
\ \ i\in I,\ \  D.v=q^{\lambda(d)}.v\}.
\end{equation*}
The subspaces $(V_q)_\lambda^+$ are defined in the obvious way. 
We say that a type 1 module is integrable if the elements 
$\bx_{i,k}^\pm$ act locally nilpotently on $V_q$ for all $i\in I$ and $k\in\bz$. 
\end{defn}

As in the classical case, one shows \cite{L2} that, if $V_q$ is integrable, then 
\begin{equation*} (V_q)_\lambda\ne 0\implies (V_q)_{\sigma(\lambda)}\ne 0\ \ 
\forall \sigma\in W.\end{equation*}
The type 1 $\bu_q^{fin}$-modules and their weight spaces are defined 
analogously. If $\lambda\in P_+$, there is a unique finite-dimensional 
irreducible $\bu_q^{fin}$-module $V_q^{fin}(\lambda)$ generated by a vector $v$ 
such that
\begin{equation*} k_i.v =q^{\lambda(h_i)}v,\ \ x_{i,0}^+.v =0,\ \ 
(x_{i,0}^-)^{\lambda(h_i)+1}.v =0,\end{equation*}
for all $i\in I$. Further, 
\begin{equation*} {\text{dim}}_{\bc(q)}(V_q^{fin}(\lambda)_\mu) 
=\dim_\bc(V^{fin}(\lambda)_\mu)\end{equation*}
for all $\mu\in P$.

From now on, we shall only consider modules of type 1.
The next result is the quantum analogue of Proposition \ref{1.2} and is proved 
in exactly the same way.
\begin{prop}{\label{qint}} Let $V_q$ be an integrable $\bu^e_q$-module. Let 
$\lambda\in P^e_+$ and assume that $0\ne v\in (V_q)^+_\lambda$. 

\begin{enumerate}
\item[(i)] $\blambda_{i, m}.v =0$ for all $i\in I$ and $|m|>\lambda(h_i)$.
\item[(ii)] $\blambda_{i,\lambda(h_i)}\blambda_{i,-m}.v 
=\blambda_{i,\lambda(h_i)-m}.v$ for all $i\in I$ and  $0\le m\le\lambda(h_i)$.
\item[(iii)] For $r\ge 1$, $s>\lambda(h_i)$, $m\in\bz$, $i\in I$, 
\begin{align*}
\left(\bX_i^-(u)^r\blambda_i^+(u)\right)_s.v  =0,\ \ &\ \
\left(\bX_{i,0}^-(u)^r\blambda_i^+(u)\right)_s.v  =0,\\
 \left(\tilde \bX^-_i(u)\blambda_i^+(u)\right)_m.v =0,\ \ &\ \ 
\left(\tilde \bH_i(u)\blambda_i^+(u)\right)_m.v =0.\ \ \ \ \ \ \ \
\end{align*}
\item[(iv)] \begin{equation*} \left(\Phi(\bX_i^-(u)^r)\blambda_i^-(u)\right)_s.v 
 =0,\ \ \ \
\left(\Phi(\bX_{i,0}^-(u)^r)\blambda_i^-(u)\right)_s.v  =0.\end{equation*}
\end{enumerate}\hfill\qedsymbol
\end{prop}

\begin{prop}{\label{fingen}}
 Let $V_q$ be an integrable type 1 $\bu_q^e$-module and assume that $\lambda\in 
P_+^e$, $0\ne v\in (V_q)_\lambda$ is such that $V_q=\bu^e_q.v$ and
 \begin{equation*}
\bx_{i,k}^+.v =0\ \forall i\in I, k\in\bz.
\end{equation*} 
Then, there exists $s_\lambda\ge 0$ such that $V_q$ is  spanned by the elements
 \begin{equation*}
(\bx_{i_1,s_1}^-)^{(l_1)}(\bx_{i_2, s_2}^-)^{(l_2)}\cdots (\bx_{i_k, 
s_k}^-)^{(l_k)}\bu_\ba(0).v\end{equation*}
for $0\le j\le k$, $l_j\ge 0$, $i_j\in I$,  $\ 0\le s_j\le s_\lambda$.   
\end{prop}

\begin{pf}
For any $N\ge 0$, let $V_N$ be the $\bc(q)$-subspace of $V_q$ spanned by the 
elements 
\begin{equation*}
(\bx_{i_1,s_1}^-)^{(l_1)}(\bx_{i_2, s_2}^-)^{(l_2)}\cdots (\bx_{i_k, 
s_k}^-)^{(l_k)}\bu_\ba(0).v
\end{equation*}
for $0\le j\le k$, $0\le s_j\le N$, $l_j\ge 0$.

By Lemma \ref{triangle}, we have $V_q =\bu_q(<)\bu_q(0).v$ and hence 
\begin{equation*} 
 V_q =\bigoplus_{\eta,m}(V_q)_{\lambda-\eta+m\delta},
\end{equation*}
where $\eta\in Q^+$, $m\in\bz$. The argument given in the proof of Proposition 
1.2(i) (but replacing the modules by their quantum analogues) shows that 
$(V_q)_{\lambda-\eta+m\delta}\ne 0$ for only finitely many 
$\eta\in Q^+$. Hence, it suffices to prove that:

\medskip\noindent {\it for each $\eta\in Q^+$, there exists $N(\eta)\ge 0$ such 
that $(V_q)_{\lambda-\eta+m\delta}\subset 
V_{N(\eta)}$ for all $m\in\bz$}.
\medskip

We proceed by induction on ${\text{ht}}\ \eta$. By Proposition \ref{qint}, we 
see 
that, if $s>\lambda(h_i)$,  $p\ge 1$,
\begin{equation*} 
\left(\bX_i^-(u)^{(p)}\blambda_i^+(u)\right)_{ps}\bu_\ba(0).v_\lambda  
=0,\end{equation*}
or equivalently that
\begin{equation}{\label{step3}}
\sum_{l=0}^{ps}\left(\bX_i^-(u)^{(p)}\right)_{l}
\blambda_i^+(u)_{ps-l}
\bu_\ba(0).v_\lambda =0.\end{equation}
If $p=1$, it follows easily from equation (\ref{step3}) by induction on $s$ that 
\begin{equation*} \bx_{i,s}^-.(V_q)_{\lambda+m\delta}\subset V_{\lambda(h_i)}
\end{equation*}
if $s>\lambda(h_i)$. To deal with the case, $s\le 0$, we apply $h_{i,s}$ to both 
sides of equation (\ref{step3}), and as in the proof of Proposition \ref{wint1} 
use the fact that $\Lambda_{i,\lambda(h_i)}$ is invertible.
Thus, the induction begins with $N(\alpha_i)=\lambda(h_i)$.  

Assume that we have proved the italicized statement above for all $\eta\in Q^+$ 
with ${\text{ht}}\ \eta<p$. We deal first with the case $\eta=p\alpha_i$, for 
some $i\in I$, $p\ge 1$. We show 
that we can take $N(p\alpha_i)=\lambda(h_i)$. For this, it suffices to prove 
that
\begin{equation} {\label{step1}}
(\bx_{i,s}^-)^{(l)}(V_q)_{\lambda+(m-sl)\delta-(p-l)\alpha_i}\subset 
V_{\lambda(h_i)} 
\end{equation}
for all $m,s\in\bz$, $p\ge l>0$. We prove this by induction on $s$ (for $p$ 
fixed), assuming first that $s\ge 0$. The induction begins 
since there is nothing to prove if $0\le s\le\lambda(h_i)$. Assume that 
$s>\lambda(h_i)$ and that the result holds for all smaller values of $s\ge 0$.

If $p>l>0$ then, by the induction on $p$, we have
\begin{equation*} 
(V_q)_{\lambda -(p-l)\alpha_i +(m-sl)\delta}\subset \sum_{s',l'} 
(\bx_{i,s'}^-)^{(l')}V_{\lambda-(p-l-l')\alpha_i+(m-sl-s'l')\delta},
\end{equation*}
where $0\le s'\le \lambda(h_i)$, $l'>0$. By Lemma \ref{order}(i),  we have 
\begin{equation*}
(\bx_{i,s}^-)^{(l)}(\bx_{i,s'}^-)^{(l')}V_{\lambda-(p-l-l')\alpha_i+
(m-sl-s'l')\delta}\subset\sum (\bx_{i,s''}^-)^{(l'')} 
V_{\lambda-(p-l'')\alpha_i+(m-l''s'')\delta},
\end{equation*}
where $s'\le s''<s$. It follows by the induction hypothesis on $s$ 
that
\begin{equation*} 
\sum (\bx_{i,s''}^-)^{(l'')} V_{\lambda-(p-l'')\alpha_i+(m-l''s'')\delta}\subset 
V_{\lambda(h_i)},\end{equation*}
which proves equation (\ref{step1}).

In the case $p=l$, we must show that 
\begin{equation*} (\bx_{i,s}^-)^{(p)}\bu_\ba(0).v \subset 
V_{\lambda(h_i)}.\end{equation*}
Since $ps>\lambda(h_i)$, by Proposition \ref{qint} we see that 
\begin{equation*}
\label{step1.5} (\bX_i^-(u)^{(p)})_{ps}\bu_\ba(0).v 
+\sum_{s'<ps}(\bX^-_i(u)^{(p)})_{s'}\blambda_{i,ps-s'}\bu_\ba(0).v 
=0.\end{equation*}
Now, by Lemma \ref{gi}(iii), we  see that for $s'<ps$ the element 
$\left(\bX_i^-(u)^{(p)}\right)_{s'}$ is a sum of terms  of the form
\begin{equation}{\label{step2}} 
(\bx_{i,s_1}^-)^{(r_1)}(\bx_{i,s_2}^-)^{(r_2)}\cdots 
\end{equation}
where $0\le s_1\le s_2\le\cdots $, $r_1,r_2,\cdots>0$, and  $s_1<s$. The 
induction hypothesis on $p$ and $s$ 
proves that $\left(\bX_i^-(u)^{(p)}\right)_{s'}\bu_\ba(0).v\subset 
V_{\lambda(h_i)}$. Finally, again by Lemma \ref{gi}(iii), we have  
\begin{equation*}
\left(\bX_i^-(u)^{(p)}\right)_{ps} = (\bx_{i,s}^-)^{(p)} + A,
\end{equation*} 
where $A$ is a sum of terms of the form in (\ref{step2}) where either $s_1<s$ or 
$s_1=s$ and $r_1<p$. If $s_1<s$ it follows as before that 
$\left(\bX_i^-(u)^{(p)}\right)_{s'}\bu_\ba(0).v\subset 
V_{\lambda(h_i)}$, and if $s_1=s$ it follows by the case $l<p$ of equation 
(\ref{step1}) proved above. Thus, the induction on $s$ 
is completed in the case $p=l$ too. 

This completes the proof of (\ref{step1}) when $s\ge 0$.
Next, consider the case when $s\le 0$. The case $p=1$ was proved above. For 
$p<l$, the same method used for $s\ge 0$ works, this time using Lemma 4.2(ii). 
Finally, for the case $p=l$, we use the relation 
\begin{equation*}\left(\Phi(\bX_i^-(u)^p)\blambda_i^-(u)\right)_{ps}.v 
=0\end{equation*}
and parts (i) and (ii) of Proposition \ref{qint}, and proceed as in the case 
$s\ge 0$. We omit the details.

 This completes the proof of the  italicized statement when $\eta$ is a multiple 
of $\alpha_i$ .

We now turn to the case of arbitrary $\eta=\sum r_i\alpha_i$ of height $p$. 
Choose $M$ so that if 
$\sum r_i<p$ then $(V_q)_{\lambda-\eta+m\delta}\subset V_M$. As in the special 
case $\eta=p\alpha_i$, to complete the induction on $p$ it suffices to prove 
that there exists $N\ge 0$ such that
\begin{equation*} 
(\bx_{i,s}^-)^{(l)}(V_q)_{\lambda-(\eta-l\alpha_i)+(m-ls)\delta}\subset 
V_{N}
\end{equation*}
for all $i\in I$, $s\in\bz$, $l>0$. Since 
${\text{ht}}(\eta-l\alpha_i)<p$, it follows that
\begin{equation*}
(V_q)_{\lambda-(\eta-l\alpha_i)+(m-ls)\delta}\subset
\sum_{j,s_j,l_j}(\bx_{j,s_j}^-)^{(l_j)}(V_q)_{\lambda-(\eta-l\alpha_i-
l_j\alpha_j)+(m-ls-l_js_j)\delta}
\end{equation*}
where $0\le s_j\le M$. Thus, we must show that there exists $N$ such that 
\begin{equation*}
(\bx_{i,s}^-)^{(l)}(\bx_{j,s_j}^-)^{(l_j)}(V_q)_{\lambda-(\eta-l\alpha_i-
l_j\alpha_j)+(m-ls-l_js_j)\delta}\subset V_N,
\end{equation*}
for all $i,j\in I$, $s,s_j\in\bz$, $0\le s_j\le M$. Assume that $s\ge 0$ (the 
case $s\le 0$ is similar). If 
$s\le M$, there is nothing to prove. Assume that we know the result for all $i$, 
and all smaller positive values of $s$. If $i=j$, then we prove exactly as in 
the case $\eta= r\alpha_i$ that we can take $N=M$. 

If $a_{ij}=0$, the result is obvious. Assume now that $a_{ij} =-1$. Then, with 
the notation in 
Proposition \ref{ij}, we see that 
\begin{equation*}
(\gamma_{s,s_j}^{i,j})^{(m)} = (-1)^m 
(\gamma_{s_j+1,s-1}^{j,i})^{(m)} =\sum_{p'=0}^m g_{p'}(\bx_{i,s-1}^-)^{(p')} 
(\bx_{j,s_{j}+1}^-)^{(m)}(\bx_{i,s-1}^-)^{(m-p')},
\end{equation*}
where the $g_{p'}\in\ba$. Using the induction hypothesis on $s$, we get that 
\begin{equation*}
(\gamma_{s,s_j}^{i,j})^{(m)}.(V_q)_{\lambda-(\eta-m\alpha_i-m\alpha_j)+m'\delta}
\subset V_{M+1}
\end{equation*} 
for all $m'\in\bz$. Now, using Proposition \ref{ij} again, we see that
\begin{equation*}
(\bx_{i,s}^-)^{(l)}(\bx_{j,s_j}^-)^{(l_j)}(V_q)_{\lambda-(\eta-l\alpha_i-
l_j\alpha_j)+(m-ls-l_js_j)\delta}\subset V_{M+1}.
\end{equation*}
This proves the result.\end{pf}

Let $I_q(0)$ be the left ideal in $\bu^e_q$ 
generated by the elements
$\blambda_{i,m}$ ($i\in I$, $|m|>\lambda(h_i)$) and 
$\blambda_{i,\lambda(h_i)}\blambda_{i,-m} 
-\blambda_{i,\lambda(h_i)-m}$ ($i\in I$, $0\le m\le\lambda(h_i)$).
Let $I_q^e(\lambda)$ be the left ideal in $\bu_q^e$ generated by 
$I_q(0)$, $\bu_q(>)$,  the elements $K_i^{\pm 1}-q^{\pm\lambda(h_i)}$ ($i\in 
I$), and $D^{\pm 1}-q^{\pm\lambda(d)}$. The ideal $I_q(\lambda)$ in 
$\bu_q(\lambda)$ is defined in the obvious way.
Let 
\begin{equation*} 
W_q(\lambda) =\bu^e_q/I_q^e(\lambda) =\bu_q/I_q(\lambda) 
\end{equation*}
be the corresponding left $\bu^e_q$-module. Let 
$w_\lambda$ be the image of $1$ in $W_q(\lambda)$. We have  
\begin{equation*} 
\bu_q(0).w_\lambda \cong\bu_q(0)/I_q(0)\end{equation*} 
as $\bu_q(0)$-modules.

\begin{prop}\label{wqint} The $\bu_q^e$-module $W_q(\lambda)$ is 
integrable.\end{prop} 
\begin{pf} In \cite{K}, Kashiwara defines an integrable $\bu_q^e$-module 
$V^{max}(\lambda)$ for all $\lambda\in P_+^e$. In fact, according to unpublished 
work of Kashiwara, 
\begin{equation*} 
V^{max}(\lambda)\cong W_q(\lambda),
\end{equation*}
and hence $W_q(\lambda)$ is an integrable $\bu^e_q$-module.  One can also give 
a direct proof that $W_q(\lambda)$ is integrable along the lines of the proof of 
Theorem \ref{wint}. One works with the presentation of $\bu^e_q$ in terms of 
Chevalley generators and the quantum version of Lemma \ref{univint} proved in 
\cite{K}, \cite{L2}. We omit the details.
\end{pf}

Given any $\bu_q^e$-module $V_q$ and a $\bu_\ba$-submodule $V_\ba$ of $V_q$ 
such that
\begin{equation*} 
V_q\cong\bc(q)\otimes_\ba V_\ba,
\end{equation*}
we set 
\begin{equation*} 
\overline {V_q}\cong \bc_1\otimes_\ba V_\ba,
\end{equation*}
where we regard $\bc_1$ as an $\ba$-module by letting $q$ act as 1. The algebra 
$\bc_1\otimes_\ba\bu_{\ba}$ is a quotient of $\bu^e$, so  
$\overline{V_q}$ is a $\bu^e$-module. Similar results hold for $\bu$-modules.

Let $\bpi_q$ be an $n$-tuple of polynomials with constant term 1 and  
coefficients in $\bc(q)$. Define $\lambda_{\bpi_q}\in\P_+^e$ and 
$\bpi_q^\pm(u)$ as in Section 2. 
Let $I_q(\bpi_q)$ be the left ideal in $\bu_q$ generated by 
$I_q(\lambda_{\bpi_q})$ and the elements
\begin{equation*}
\left(\blambda_i^\pm(u)-\pi_i^\pm(u)\right)_s\ \ (i\in I, s\ge 0).
\end{equation*} 
Set $W_q(\bpi_q) = \bu_q/I_q(\bpi_q)$.

\begin{lem} $W_q(\bpi_q)$ is a finite-dimensional $\bu_q$-module.\end{lem}
\begin{pf}  This is proved in the same way as the corresponding result for 
$\bu$. We use  Proposition \ref{fingen} instead of Proposition \ref{wint1}.
\end{pf}

One has the following analogue of Proposition \ref{universal} for the modules 
$W_q(\lambda)$ and $W_q(\bpi_q)$. We omit the proof, which is entirely similar 
to that of Proposition \ref{universal}.

\begin{prop} 
\begin{enumerate}
\item[(i)] Let $V_q$ be any integrable $\bu_q^e$-module generated by an element 
of $(V_q)_\lambda^+$. Then, $V_q$ is a quotient of $W_q(\lambda)$.
\item[(ii)] Let $V_q$ be a finite-dimensional quotient $\bu_q$-module 
of $W_q(\lambda)$, and assume that ${\text{dim}}\ (V_q)_\lambda =1$. Then, $V_q$ 
is a quotient of $W_q(\bpi_q)$ for some choice of $\bpi_q$.
\item[(iii)] Let $V_q$ be  finite-dimensional $\bu_q$-module generated by an 
element of $(V_q)_\lambda^+$ and such that ${\text{dim}}\ (V_q)_\lambda =1$.   
Then, $V_q$ is a quotient of $W_q(\bpi_q)$ for some  $\bpi_q$.
\end{enumerate}
\end{prop}  

\begin{defn} We call $\bpi_q$ integral if the polynomials
$\pi_i^\pm(u)$ have coefficients in $\ba$ for all $i\in I$. Equivalently, for 
all $i\in I$, $\pi_i(u)$ has coefficients in $\ba$ and the coefficient of the 
highest power of $u$ should lie in $\bc^\times q^\bz$. Let $\overline{\bpi_q}$ 
be the 
$n$-tuple of polynomials with coefficients in $\bc$ and constant term one 
obtained from $\bpi_q$ by evaluating its coefficients at $q=1$. 
\end{defn}

For any $\bpi_q$, $W_\ba(\bpi_q) =\bu_\ba.w_{\bpi_q}$ is a $\bu_\ba$-module.  

\begin{lem} Assume that $\bpi_q$ is integral.

\begin{enumerate}
\item[(i)]  $W_\ba(\bpi_q)$ is a free $\ba$-module and we have
\begin{equation*}
W_q(\bpi_q)\cong \bc(q)\otimes_\ba W_\ba(\bpi_q).
\end{equation*}
\item[(ii)]  The $\bu$-module $\overline{W_q(\bpi_q)}$ is a quotient of 
$W(\overline{\bpi_q})$.
\end{enumerate}

\end{lem}
\begin{pf} Since $W_\ba(\bpi_q)$ is a quotient of $W_q(\lambda_{\bpi_q})$,
it follows from Proposition \ref{fingen} that $W_\ba(\bpi_q)$ is a 
finitely-generated $\ba$-module. Since it is clearly a torsion-free 
$\ba$-module, part (i) follows. 

To prove (ii), observe that the defining relations of $W_q(\bpi_q)$ specialize 
to those of $W(\overline\bpi)$. The 
result now follows from Proposition \ref{univint}. 
\end{pf}

The $\bu_q$-module $W_q(\bpi_q)$ has a unique irreducible quotient 
$V_q(\bpi_q)$. Let $v_{\bpi_q}$ be the image of $w_{\bpi_q}$ and set
\begin{equation*} V_\ba(\bpi_q) =\bu_\ba.v_{\bpi_q}.\end{equation*}
If $\bpi_q$ is integral, $V_\ba(\bpi_q)$ is a $\bu_\ba$-module and is free as an 
$\ba$-module
(since it is torsion-free and the quotient of a finitely-generated 
$\ba$-module).
Let $V(\bpi)$ be the unique irreducible quotient of the $\bu$-module $W(\bpi)$. 

\begin{lem} The $\bu$-module $V(\overline{\bpi_q})$ is a quotient of  
$\overline{V_q(\bpi_q)}$.\end{lem}
\begin{pf} The module $\overline{V_q({\bpi_q})}$ has an unique irreducible 
quotient $V$. By Lemma 4.6(ii), $V$ is a quotient of $W(\overline{\bpi_q})$ and 
hence by 
uniqueness $V\cong  V(\overline{\bpi_q})$. 
\end{pf}

We have thus proved the following statement. Assume that $\bpi_q$ is integral. 
Then, we have a commutative diagram of surjective $\bu$-module homomorphisms
 \begin{eqnarray*}
 W(\overline{\bpi_q})&\longrightarrow&\overline{W_q(\bpi_q)}\\
 \downarrow&&\downarrow\\
 V(\overline{\bpi_q})&\longleftarrow&\overline{V_q(\bpi_q)}.
 \end{eqnarray*}

\vskip9pt\noindent{\bf Conjecture.} If $\bpi_q=\bpi$ has coefficients in $\bc$, 
the natural map 
$W(\bpi)\to\overline{V_q(\bpi)}$ is an isomorphism of $\bu$-modules, and hence 
$W_q(\bpi)\cong V_q(\bpi)$. \hfill\qedsymbol
\vskip9pt
In Section 6, we prove this conjecture when $\frak g=sl_2$.

\section{An irreducibility criterion.} 

In this section we return to the classical case and obtain a criterion for the 
irreducibility of the modules $W(\bpi)$. 

 For any $a\in\bc^{\times}$, and any $\ug$-module $V$, define a $\bu$-module 
structure on $V$ by
\begin{equation*} (x\otimes t^r). v= a^r x.v\end{equation*}
for $x\in\frak g$, $r\in\bz$, $v\in V$. Let $V(a)$ denote the corresponding 
$\bu$-module.

For $i\in I$, $a\in\bc^\times$, we denote by $W(i,a)$ the $\bu$-module 
corresponding to the $n$-tuple $\bpi$ of polynomials defined by
\begin{equation*}
\pi_j(u) =1 \ {\text{if}}\ j\ne i,\ \ \pi_i(u) =1-au,
\end{equation*}
and denote $w_\bpi$ by $w_{i,a}$. Clearly, $V^{fin}(\omega_i)(a)$ is the 
irreducible quotient of $W(i,a)$.  We set $V^{fin}(\omega_i)(a) = V(i,a)$.

For $i\in I$ and $a\in\bc(q)^{\times}$, the $\bu_q$-modules $W_q(i,a)$ and 
$V_q(i,a)$ are defined similarly.

We need the following result, due to \cite{CP3} for $\frak g$ of type 
$sl_2$ and due to 
\cite{K2} and \cite{FM} in general.
\begin{prop}{\label{irrtensor}} Let $r\ge 1$, $a_1,\cdots 
,a_r\in\bc(q)^{\times}$, $i_1,i_2,\cdots,
i_r\in I$. There is a finite set $S\subset\bc(q)^\times$ (depending on 
$i_1,\cdots,i_r$) such that the tensor product 
\begin{equation*} 
V_q(i_1,a_1)\otimes V_q(i_2,a_2)\otimes\cdots\otimes V_q(i_r, a_r),
\end{equation*} 
is irreducible if $a_l/a_m\notin S$ for all $l,m=1,2,\ldots,r$. If $\frak 
g=sl_2$, $S=\{q^{\pm 2}\}$.\hfill\qedsymbol\end{prop}

\begin{prop} {\label{irrfund}}For $i\in I$, $a\in \bc^\times$, $W(i,a)\cong 
V(i,a)$ 
if and only if  $r_i=1$.  \end{prop}
 \begin{pf} The proof rests on the following fact, which can be established by a 
case by case check: $r_i=1$ if and only if there exists $\mu\in P_+$ with 
$0\ne\omega_i-\mu\in Q_+$ such that $V^{fin}(\mu)$ occurs as a component of 
$V^{fin}(\theta)\otimes V^{fin}(\omega_i)$.

Suppose first that $r_i>1$. Let $\mu\in P_+$ have the above property. For 
$x\in\frak g$, $m\in\bz$, $v\in V^{fin}(\omega_i)$, $v'\in V^{fin}(\mu)$, define
\begin{equation*}
x_m.(v,v')=a^m(x.v, \ \  m\, pr(x\otimes v)+x.v'),
\end{equation*}
where $pr:V^{fin}(\theta)\otimes V^{fin}(\omega_i)\to V^{fin}(\mu)$ is the 
$\frak g$-module projection. It is straightforward to check that this defines a 
$\bu$-module structure on $V^{fin}(\omega_i)\oplus V^{fin}(\mu)$, and that this 
$\bu$-module is generated by the highest weight vector in $V^{fin}(\omega_i)$. 
It is therefore a quotient of $W(i,a)$.

 To prove the converse, notice that, as a $\frak g$-module, $W(i,a)$ is 
completely reducible, and hence
\begin{equation*} W(i,a) \cong V^{fin}(\omega_i)\oplus\bigoplus_{\mu<\omega_i} 
V^{fin}(\mu)^{m_\mu},\end{equation*}
 where $m_\mu$ is the multiplicity with which $V^{fin}(\mu)$ occurs in $W(i,a)$. 
 Consider the map $\loopg\otimes W(i,a)\to W(i,a)$  given by 
\begin{equation*}
x_n\otimes v\mapsto x_n.v.
\end{equation*}
 This is clearly a map of $\frak g$-modules, where we regard $\loopg$ as a 
module for $\frak g$ through the adjoint representation. For each $m\in\bz$, 
consider the restriction of this map to $(\frak g\otimes t^m)\otimes 
V(\omega_i)$. Since $\frak g\otimes t^m\cong V^{fin}(\theta)$ as $\frak 
g$-modules, 
we have a $\frak  g$-module  map $V^{fin}(\theta)\otimes V(i,a)\to W(i,a)$. By 
the fact stated above, this map takes $(\frak g\otimes 
t^m)\otimes V(\omega_i)$ into the $\bu(\frak g)$-submodule $V^{fin}(\omega_i)$ 
of 
$W(i,a)$ for all $m\in\bz$. This proves that $V^{fin}(\omega_i)$ is a 
$\bu$-submodule of $W(i,a)$ and 
hence (since $w_{i,a}\in V^{fin}(\omega_i)$) is equal to $W(i,a)$.
\end{pf}

\begin{remark} 
The same criterion $r_i=1$ occurs, for the same reason, in Drinfeld's work 
on finite-dimensional representations of Yangians, \cite{Dr1}. See also 
\cite[Proposition 12.1.17]{CP4}.
\end{remark}
 
We can now state the main result of this section.

\begin{thm}{\label{irred}} 
Let $\bpi=(\pi_1,\ldots,\pi_n)$ be an $n$-tuple of polynomials in $\bc[u]$ with 
constant coefficient one. Then, the $\bu$-module $W(\bpi)$ is irreducible if 
and only if $\pi_\theta$ has distinct roots.
\end{thm}
\begin{pf} Assume that $\pi_\theta$ has distinct roots. By Lemma \ref{pibeta1}, 
it follows that $\pi_i=1$ if $r_i\ne 1$. Let 
\begin{equation*}
I'=\{i\in I: r_i=1\}.
\end{equation*} 
If $i\in I'$ then $\pi_i$ must have distinct roots, and for any $i,j\in I'$, 
$i\ne j$, the 
polynomials  $\pi_i$ and $\pi_j$ 
must be relatively prime. Hence, by Theorem \ref{wtensor} and Proposition 
\ref{irrfund}, it follows that
\begin{equation*} 
W(\bpi)\cong\bigotimes_{i\in I',a_{ij}\in\bc^\times}
W(i,a_{ij})\cong\bigotimes_{i\in I',a_{ij}\in\bc^\times}V(i,a_{ij}),
\end{equation*} 
where the $a_{ij}$ are all distinct. 
By Proposition \ref{vtensor}, we see  that the second  tensor product in the 
preceding 
equation is an irreducible
$\loopg$-module.

For the converse, suppose that $\pi_\theta$ has repeated roots. By Theorem 
\ref{wtensor}, it follows that $W(\bpi)$ is isomorphic to a tensor product
of modules $W(\bpi^a)$, where $\bpi^a$ is an $n$-tuple of polynomials 
such that $\pi^a_\theta =(1-au)^{m}$ for some $a\in\bc^\times$ and $m\ge 1$, and 
where $m>1$ for at least one value of $a$.  Thus, it suffices to prove the 
theorem in the case where   $\pi_\theta(u)=(1-au)^m$ with $a\in\bc^\times$ and 
$m>1$. From now on, we shall assume that we
 are in this case. 

To prove that $W(\bpi)$ is not isomorphic to $V(\bpi)$ as $\loopg$-modules, 
recall that by Proposition \ref{vtensor}, we have $V(\bpi)\cong 
V^{fin}(\lambda_\bpi)$ as $\frak g$-modules.  Hence, it suffices to prove that 
$W(\bpi)$ is reducible 
as a $\frak g$-module.  

Let $\bpi_q$ be an $n$-tuple of polynomials with constant term 1 such that
\begin{equation*}
\pi_i(u) =(1-a_{i,1}u)(1-a_{i,2}u)\cdots(1-a_{i,m_i}u),
\end{equation*}
where $a_{i,j} =aq^{l_{ij}}$, for  some $l_{ij}\in\bz$. Let $v_{i,j}$ be the 
highest weight vector in  $V_q(i,a_j)$ and consider 
\begin{equation*}
\bV =V_q(1,a_{1,1})\otimes 
V_q(1,a_{1,2})\cdots \otimes V(1,a_{1, m_1})\otimes\cdots \otimes 
V_q(n,a_{n,1})\otimes \cdots\otimes V_q(n,a_{n,m_n}).
\end{equation*}
Let $\bv=v_{1,1}\otimes v_{1,2}\otimes\cdots\otimes v_{n,m_n}$ and 
set
\begin{equation*} 
Z_q(\bpi_q) =\bu_q.\bv\subset\bV,\ \ \ Z_\ba(\bpi_q) =\bu_\ba.\bv.
\end{equation*} 
Since $Z_q(\bpi_q)$ is a quotient of $W_q(\bpi_q)$,  and $\bpi_q$ is integral, 
it follows that 
\begin{equation*}
Z_q(\bpi_q)\cong Z_\ba(\bpi_q)\otimes_{\ba}\bc(q),
\end{equation*}
so we can define the $\bu$-module $\overline{Z_q(\bpi_q)}= 
Z_\ba(\bpi_q)\otimes_\ba\bc_1$.
Clearly, $\overline{Z_q(\bpi_q)}$ is a quotient of $W(\bpi)$ and hence it 
suffices 
to show that $\overline{Z_q(\bpi_q)}$ is reducible as a $\frak g$-module.  

Suppose first that  $m_{i_0}\ge 2$ for some ${i_0}\in I$. Take $l_{ij}= 0$ for 
all  
$i\in I$ and $j=1,\cdots ,m_{i}$. Let $\bu_q^{i_0}$ be the subalgebra of $\bu_q$ 
generated by $K_{i_0}^{\pm 1}$ and $x_{{i_0},k}^\pm$ for $k\in\bz$. Consider the 
$\bu_q^{i_0}$-module 
\begin{equation*}
Z_q^{i_0}(\bpi_q) =\bu_q^{i_0}.\bv.
\end{equation*}
Let $\omega$ be the fundamental weight of $sl_2$. Then, by Proposition 
\ref{irrtensor}, we know 
that $V_q(\omega, a)^{\otimes m_{i_0}}$ is irreducible and hence is  a quotient 
of 
$Z_q^{i_0}(\bpi_q)$. Clearly,
\begin{equation*} 
\dim(Z_q(\bpi_q)_{\lambda-\alpha_{i_0}})\ge\dim(V(\omega, a)^{\otimes 
m_{i_0}})_{m_{i_0}\omega -\alpha} = m_{i_0}
\end{equation*}
hence
\begin{equation*}
\dim(\overline{Z_q(\bpi_q)}_{\lambda-\alpha_{i_0}}\ge m_{i_0}>1.
\end{equation*}
On the other hand, $V(\bpi)$ is a quotient $\bu$-module of 
$\overline{Z_q(\bpi_q)}$, since $\bpi_q=\bpi$, and $V(\bpi)\cong 
V^{fin}(\lambda_\bpi)$ as $\bu^{fin}$-modules from above. But
\begin{equation*}
\dim(V^{fin}(\lambda_\bpi)_{\lambda_\bpi-\alpha_{i_0}})=1.
\end{equation*}
Hence, $\overline{Z_q(\bpi_q)}$ is reducible as a $\bu^{fin}$-module.

We can therefore assume that each $m_i=0 \ {\text{or}} \ 1$, and that at least 
one $m_i=1$.  Consider first the case $m_{i_0} =m_{i_1} =1$ with 
$i_0< i_1$, and $m_j=0$ for all $i_0<j<i_1$. Set $J=\{i_0,i_0+1,\ldots,i_1\}$. 
By a suitable choice of the numbering, we can assume that the corresponding 
diagram subalgebra $\frak g^J$ of $\frak g$ is of type $A_{|J|}$. Let $\bu_q^J$ 
be the subalgebra of $\bu_q$ generated by 
$K_i^{\pm 1}$ and $x_{i,k}^\pm$ for $k\in\bz$ and $i\in J$. Define 
\begin{equation*}{Z_q^J}(\bpi_q)={\bu_q^J}.\bv.\end{equation*}
By Proposition 5.2, the $\bu_q^J$-module $V_q^J(i_0, a_{i_0})\otimes V_q^J(i_1, 
a_{i_1})$ is 
irreducible except for  finitely many values of the ratio $a_{i_0}/a_{i_1}$. 
Since each $a_i$ can be chosen from the infinite set  $\{aq^m:m\in\bz\}$, we can 
assume that $a_{i_0}$ and $a_{i_1}$ are chosen so that $V_q^J(i_0, 
a_{i_0})\otimes 
V_q^J(i_1, a_{i_1})$ is an irreducible $\bu_q^J$-module and hence a quotient of 
$Z_q^J(\bpi_q)$. If $\theta_J$ is the highest root of the subdiagram $J$, then
\begin{align*}  
{\text {dim}} (Z_q(\bpi_q)_{\lambda_{\bpi_q}-\theta_J}) &\ge 
{\text {dim}}(Z^J_q(\bpi_q)_{\lambda_{\bpi_q}-\theta_J})\\& \ge {\text 
{dim}}((V_q^J(i_0, a_{i_0})\otimes V_q^J(i_1, 
a_{i_1}))_{\lambda_{\bpi_q}-\theta_J})\\& 
\ge\dim((V_q^{fin,J}(\omega_{i_0})\otimes 
V_q^{fin,J}(\omega_{i_1}))_0)=|J|+1\end{align*}
(in an obvious notation). On the other hand, $V(\bpi)$ is a quotient 
$\bu$-module of $\overline{Z_q(\bpi_q)}$ and
\begin{equation*}
V(\bpi)_{\lambda_\bpi-\theta_J}=V^J(\bpi)_{\lambda_\bpi-\theta_J}=V^{fin,J}
(\lambda_\bpi)_{\lambda_\bpi-\theta_J},
\end{equation*}
which has dimension $|J|$. Hence, $V(\bpi)$ is a proper quotient of 
$\overline{Z_q(\bpi_q)}$. This shows that $Z_q(\bpi_q)$ is not isomorphic to 
$V_q(\lambda_{\bpi_q})$ and 
hence is reducible as a $\bu_q(\frak g)$-module.

It remains to consider the case when exactly one $m_i=1$, say $m_{i_0} =1$ and 
all other $m_i=0$. Since $\pi_\theta$ has repeated roots, this means that 
$r_{i_0}\ne 1$. By Proposition \ref{irrfund} we know that $W(\bpi)=W(i,a)$ is  
reducible. 

This completes the proof of the theorem.\end{pf}

\section{ The $sl_2$ case.}
In this section, $\frak g =sl_2$. Let $\omega$ be the fundamental weight, 
$\alpha$ the positive root, and set $x^\pm=x_{\pm\alpha}$, $h=h_\alpha$. Let 
$\pi$ be a polynomial with coefficients in $\bc$ and constant term 1. When $\pi 
=1-au$, denote $V_q(\pi)$ by $V_q(a)$, and define $V(a)$ similarly. Note that 
these modules are two-dimensional over $\bc(q)$ and $\bc$, respectively.

Set $V=V^{fin}(\omega)$ and let $L(V) =V\otimes \bc[t,t^{-1}]$ be the obvious 
$\bu^e$-module, given by
\begin{equation*} 
x_r.(v\otimes t^s) =x.v\otimes t^{r+s},\ \ d.(v\otimes t^r) = 
rv\otimes t^r,
\end{equation*}
for $r,s\in\bz$, $x\in\frak g$ and $v\in V$.

Set $\cal{P}_m = \bc[t_1^{\pm 1},t_2^{\pm 1},\cdots ,t_m^{\pm 1}]$. Let 
$\Sigma_m$ be the symmetric group on $m$ letters and let $\cal{P}^{\Sigma_m}$ be 
the subalgebra of polynomials invariant under the obvious action of $\Sigma_m$.

Let $S^m(L(V))$ be the symmetric part of the $m$--fold tensor product 
$T^m(L(V))$ of $L(V)$. Then, $T^m(L(V))$ is a $\bu^e$--module in the obvious 
way, and $S^m(L(V))$ is a $\bu^e$--submodule. Moreover, as vector spaces,
\begin{equation*} T^m(L(V))\cong V^{\otimes m}\otimes \cal{P}_m,\end{equation*} 
and so $T^m(L(V))$ is a right module for $\cal{P}$ by right multiplication.
This induces a right $\cal{P}_m^{\Sigma_m}$-action on $S^m(L(V))$.

The left $\bu^e$--module $W(m\omega)$ is also a right 
$\cal{P}^{\Sigma_m}$--module. In fact, by equation (\ref{rightsym}), 
$W(m\omega)$ is  a right $\bu(0)/I_{m\omega}(0)$--module, i.e., a right module 
for the algebra $\bc[\Lambda_1,\cdots, \Lambda_m,\Lambda_m^{-1}]$. But this 
algebra is isomorphic to $\cal{P}^{\Sigma_m}$ by taking $\Lambda_r$ to the 
$r^{th}$ elementary symmetric function of $t_1,\cdots ,t_m$.

 In this section, we shall prove the following two theorems.
\begin{thm}{\label{rankw}} 
As left $\bu^e$-modules and as right $\cal{P}_m^{\Sigma_m}$-modules, we have
\begin{equation*} 
W(m\omega) \cong S^m(L(V)).
\end{equation*} 
\end{thm}
To prove Theorem \ref{rankw} we shall need
\begin{thm}{\label{dimw}} Let $\pi(u)$ be a polynomial with coefficients in 
$\bc$ and constant term 1. Then,
the dimension of $W(\pi)$ is $2^{{\text{deg}}\pi}$. In fact, 
\begin{equation*} 
W(\pi)\cong \overline{V_q(\pi)},
\end{equation*}
as $\bu$-modules, and
\begin{equation*} 
V_q(\pi)\cong\bigotimes_a V_q(a)
\end{equation*}
where $a^{-1}$ runs over the set of roots of $\pi$ counted with multiplicity.  
\end{thm}
It {\it{does not}} follow from this result that $W(\pi)\cong\bigotimes_a 
V(a)$. In fact, this is false except when ${\text{deg}}\pi=1$. The point is that 
the 
$\ba$-form of $\bigotimes_a V_q(a)$ is not the tensor product of the 
$\ba$-forms of $V_q( a)$ (in fact, the former is  a proper subset of the 
latter, unless ${\text{deg}}\pi =1$).
We note the following corollary.\begin{cor} For any $\pi(u)$ as in Theorem 
\ref{dimw}, 
we have 
 $W_q(\pi)\cong V_q(\pi)$ as $\bu_q$--modules.\end{cor}
\begin{pf} Since $V_q(\pi)$ is a quotient of $W_q(\pi)$ it suffices to prove 
that \begin{equation*} {\text{dim}}_{\bc(q)}W_q(\pi)\le 
2^{{\text{deg}}\pi}.\end{equation*}
But this is clear from Theorem \ref{dimw},  since $\overline{W_q(\pi)}$ is a 
quotient of $W(\pi)$, so \begin{equation*} {\text{dim}}_{\bc(q)}W_q(\pi) 
={\text{dim}}_\bc\overline{W_q(\pi)}\le {\text{dim}}_\bc W(\pi) 
=2^{{\text{deg}}\pi}.\end{equation*}
\end{pf}

Assume Theorem \ref{dimw} for the moment.
To prove Theorem \ref{rankw}, we begin with the following trivial lemma.

\begin{lem}{\label{sint}} 
The $\bu^e$-module $S^m(L(V))$ is integrable. \hfill\qedsymbol
\end{lem}

Let $\{v_+,v_-\}$ be the usual basis of $V$, so that
\begin{equation*} 
x^\pm.v_\pm =0,\ \ x^\pm.v_\mp = v_\pm,\ \ hv_\pm =\pm v_\pm. 
\end{equation*}
For $0\le r\le m$, $l_1,\cdots ,l_m\in\bz$, define
\begin{equation*} 
\bv_{(l_1,l_2,\cdots,l_r), (l_{r+1},\cdots,l_m)} = 
\sum_{\sigma\in \Sigma_{m}}\ v_{\sigma(1)}t^{l_{\sigma(1)}}\otimes\cdots\otimes 
v_{\sigma(m)}t^{l_{\sigma(m)}},
\end{equation*}
where we set 
\begin{align*} v_s& =v_-, \ \ \text{if}\ 1\le s\le r,\\ v_s &=v_+,\ \ \text{if}\ 
r+1\le s\le m.\end{align*}
Clearly the set\begin{equation*}\{\bv_{(l_1,l_2,\cdots ,l_r), (l_{r+1},\cdots, 
l_m)}: 0\le r\le m, \ l_1,\cdots ,l_m \in\bz\},\end{equation*} is a  
$\bc$--basis of $S^m(L(V))$.

\begin{lem}{\label{free}} The  $\cal{P}_m^{\Sigma_m} $--module $S^m(L(V))$    is 
free of rank $2^m$.\end{lem}
\begin{pf} For $0\le r\le m$, let $S^m(L(V))_{r}$ be the subspace of $S^m(L(V))$ 
spanned by the elements
\begin{equation*}\{\bv_{(l_1,l_2,\cdots ,l_r), (l_{r+1},\cdots, l_m)}: \ 
l_1,\cdots ,l_m \in\bz\}.\end{equation*}
Clearly, $S^m(L(V))_r$ is a right $\cal{P}_m^{\Sigma_m}$--submodule of 
$S^m(L(V))$. It is easy to see that it is  isomorphic to the  
$\cal{P}_m^{\Sigma_m}$--module $\cal{P}_m^{\Sigma_r\times \Sigma_{m-r}}$, 
consisting of the elements in $\cal{P}_m$ invariant under permutation of the 
first $r$ and the last $m-r$ variables.. But it is well known that the latter 
module is free of rank $\left(\begin{matrix} m\\r\end{matrix}\right)$. This 
proves the lemma.
\end{pf}

\begin{lem}{\label{quot}} The assignment $w_{m\omega}\to \bv_+^{\otimes m}$ 
extends  to a well defined surjective homomorphism $W(m\omega)\to S^m(L(V))$ of 
left $\bu^e$--modules and right $\cal{P}_m^{\Sigma_m}$--modules.\end{lem} 

\begin{pf}
It follows by Proposition 
\ref{universal}(i)  that there exists a $\bu^e$-module homomorphism $\phi: 
W(m\omega)\to S^m(L(V))$ that takes $w_{m\omega}$ to $\bv=v_+^{\otimes m}$. It 
is trivial to check that $\phi$ is  also a map of right   
$\cal{P}_m^{\Sigma_m}$--modules.
To show that $\phi$ is surjective, it is enough to prove that 
\begin{equation} 
{\label{sym}}S^m(L(V)) =\bu^e.\bv.
\end{equation} 
We prove by induction on $r$ that 
\begin{equation}{\label{6}} 
\bv_{(l_1,l_2,\cdots,l_r), (l_{r+1},\cdots,l_m)}\subset \bu^e.\bv
\end{equation}
for all $l_1,\cdots ,l_m\in\bz$. 
Consider the case $r=0$.  For any 
$k_1,k_2,\cdots,k_m\in\bz$, we have  
\begin{equation*} 
(x^+_0)^mx^-_{k_1}x^-_{k_2}\cdots x^-_{k_m}. \bv  
=\sum_{\sigma\in\Sigma_m} v_+t^{k_{\sigma(1)}}\otimes 
v_+t^{k_{\sigma(2)}}\otimes\cdots\otimes v_+t^{k_{\sigma(m)}},\end{equation*}
which proves (\ref{6}) in this case. The case $r=m$ can be done similarly, since
the element $\bv^- =  v_-^{\otimes m} ={\frac{1}{ m!}} (x^-)^m.\bv\in \bu^e$. 

Assuming the result for $r$, we prove it for $r+1$. For this we shall proceed by an  induction on
\begin{equation*} N=\#\{j: j\ge r+1, \ \ l_j\ne 0\}.\end{equation*}
Now,   
\begin{equation*} 
x_k^-.\bv_{(l_1,l_2,\cdots,l_r), (l_{r+1},\cdots, 
l_m)}=\sum_{s=r+1}^m \bv_{(l_1,l_2,\cdots,l_r,l_s+k), (l_{r+1},\cdots, 
\hat{l}_s,\cdots,l_m)}.
\end{equation*}
Taking $l_s=0$ for all $s>r$, we get that 
\begin{equation*}
\bv_{(l_1,l_2,\cdots,l_r,k), (0,\cdots,0)}\in\bu^e.\bv,
\end{equation*}
for all $k\in\bz$, proving our assertion when $N=0$. Assume the result for 
$N-1$. We have to show that 
\begin{equation*}
\bv_{(l_1,l_2,\cdots,l_r,l_{r+1}),
(l_{r+2},\cdots,l_{r+N+1},0,\cdots,0)}\in\bu^e.\bv.
\end{equation*}
Now  
\begin{equation*}
x_{l_{r+1}}^-.\bv_{(l_1,l_2,\cdots, l_r), (l_{r+2},\cdots, 
l_{r+N+1}, 0,\cdots ,0 )}\end{equation*} is in $\bu^e.\bv$ by the induction 
hypothesis on $r$ , and is 
a sum of the term
\begin{equation*}
(m-r)\bv_{(l_1,\cdots ,l_{r+1}),(l_{r+2},\cdots,l_{r+N+1},0, \cdots ,0)}
\end{equation*} 
and terms of type
\begin{equation*}
\bv_{(l_1,\cdots,l_r,l_{s}+l_{r+1}),(l_{r+2},\cdots 
\hat{l}_s,\cdots, l_{r+N+1},0,\cdots,0)},
\end{equation*} 
for $r+2\le s\le r+N+1$. Since, by the  induction hypothesis on $N$, all the 
terms of the 
second type are in $\bu^e.\bv$, it follows that 
\begin{equation*}
\bv_{(l_1,\cdots ,l_{r+1}),(l_{r+2},\cdots, l_{r+N+1},0,\cdots,0)}\in\bu^e.\bv.
\end{equation*}  
This completes the proof that $\phi$ is surjective.\end{pf} 

\noindent {\it{Proof of Theorem \ref{rankw}}}. 
Let $K$ be the kernel of the homomorphism $W(m\omega)\to S^m(L(V))$ given by 
Lemma \ref{quot}. Since $S^m(L(V))$ is a free, hence projective, right 
$\cal{P}_m^{\Sigma_m}$--module by Lemma \ref{free}, it follows that
\begin{equation*} W(m\omega)=S^m(L(V))\oplus K,\end{equation*}
as right $\cal{P}_m^{\Sigma_m}$--modules.

Let $\frak m$ be any maximal ideal in $\cal{P}_m^{\Sigma_m}$. Identifying 
$\cal{P}_m^{\Sigma_m}$ with $\bu(0)/I_{m\omega}(0)$ as described earlier in this 
section, it is clear that
\begin{equation*} \frak m=I_\pi(0)/I_{m\omega}(0),\end{equation*}
for some polynomial $\pi$ with constant term 1 such that ${\text{deg}}\pi =m$. 
It follows that
\begin{equation*} W(m\omega)/W(m\omega)\frak m\cong W(\pi)\end{equation*}
as vector spaces over $\bc$, and hence has dimension $2^m$. On the other hand, 
by Lemmma \ref{free}, $S^m(L(V))/S^m(L(V))\frak m$ also has dimension $2^m$ over 
$\bc$. It follows that $K/K\frak m =0$. Since this holds for all maximal ideals 
$\frak m$, Nakayama's lemma implies that $K=0$, proving the theorem.

\hfill\qedsymbol
\medskip

The rest of the section is devoted to proving Theorem  \ref{dimw}. First, 
observe that, in view of Theorem \ref{wtensor}, it suffices to consider the case 
when $\bpi(u) = (1-au)^m$ for some $a\in\bc^\times$. Since we have a surjective 
map $W(\bpi)\to\overline {V_q(\bpi)}$, it suffices to prove that 
\begin{equation}{\label{dimless}} 
{\text{dim}}_\bc\ W(\bpi)\le{\text{dim}}_\bc\overline{V_q(\pi)}=2^m.
\end{equation}

For $a\in\bc^\times$, let $\tau_a: \frak g\otimes\bc[t]\to\frak g\otimes\bc[t]$ 
be the Lie algebra automorphism obtained by extending the assignment
\begin{equation*} 
x\otimes t^k\to x\otimes (t-a)^k,\ \ \forall \ x\in\frak g,\ k\ge 0.
\end{equation*} 
Set
\begin{equation*} 
X_a^-(u) =\tau_a(X_{\alpha,0}^-(u)),\ \ 
\Lambda_a^+(u)=\tau_a(\Lambda_\alpha^+(u)) 
={\text{exp}}\left(-\sum_{k=1}^\infty\frac{h\otimes(t-a)^k}{k}u^k\right).
\end{equation*}
It is easy to see, using the relation between the $\Lambda_{m}$ and $h_m$, that
\begin{equation}
{\label {hk}} h_k.w_\bpi=ma^k w_\bpi,
\end{equation}
or equivalently that 
\begin{equation*} h\otimes (t-a)^k . w_\bpi =0,\end{equation*} for all $k\ge 0$.
It follows that
\begin{equation}{\label{lwpi0}} 
(\Lambda_a^+)_k.w_\bpi =0,\ \ \forall\ k> 0.
\end{equation}
Further, using Lemma \ref{gar}(ii) and observing that the identity there is 
actually 
an identity in $\bu(\frak g\otimes \bc [t])$, we get by applying $\tau_a$ that 
\begin{equation*} \tau_a(x_1^+)^{(r)}(x_0^-)^{(s)} = 
(-1)^r\left(X_a^-(u)^{(s-r)}\Lambda_a^+(u)\right)_s 
\mod\bu\bu(>)_+,
\end{equation*} for $s\ge r\ge 1$.
Together with (\ref{lwpi0}), it follows that \begin{equation}{\label {xpi0}}
\left(X_a^-(u)^{(s-r)}\right)_s.w_\bpi =0 \ \ \ 
\forall\ \  r\ge 1, \  s\ge m+1.
\end{equation}
In particular, this means that \begin{equation*}(x^-\otimes (t-a)^{s}).w_\bpi 
=0\ \ \forall\ \  s\ge m.\end{equation*} 
Let $\bu^+_a(<)$ be the commutative subalgebra of $\bu$ generated by the 
elements $\tau_a(x^-_k)$ for all $k\ge 0$, and let $I_a(m)$ be the ideal 
in $\bu_a^+(<)$ generated by the elements $\left(X_a^-(u)^{(s-r)}\right)_s$ for 
all $r\ge 1$, 
$s\ge m+1$. 
\begin{lem}{\label{6.4}} The assignment  $u\to u.w_\bpi$ induces a surjective 
map of vector 
spaces
$\bu^+_a(<)/I_a(m)\to W(\bpi)$.
\end{lem}
\begin{pf} The map is well defined by (\ref{xpi0}). It is obviously surjective 
because the polynomials $(t-a)^k$ for $k\ge 0$ are a basis of $\bc[t]$.\end{pf}

Thus, to prove (\ref{dimless}) it suffices to show that the dimension of 
$\bu^+_a(<)/I_a(m)$ is at most $2^m$. It is convenient to reformulate the 
problem as follows. 

Let $R_m =\bc[z_0,\cdots ,z_{m-1}]$ be the polynomial algebra in $m$ variables. 
For 
$0\le j<m$, set
\begin{equation*} 
Z_j(u)=\sum_{i=j}^{m-1}z_iu^{i-j+1}\in R_m[u].
\end{equation*}
Let $J_m$ be the ideal in $R_m$ generated by the elements 
$\left(Z_0(u)^r\right)_s$, for 
$r\ge 1$, $s\ge m+1$. It is trivial to see that
\begin{equation*} R_m/J_m\cong \bu_a^+(<)/I_a(m),\end{equation*}
via the map $\tau(x_k^-)\to z_k$.
It is clear that (\ref{dimless}) is now a consequence of the following 
proposition and Lemma \ref{6.4}.

\begin{prop}{\label{basis}} For $m\ge r> 0$, let
\begin{equation*}
\cal{B}_{m,r} =\{z_{i_1}z_{i_2}\cdots z_{i_r}: 0\le i_1\le 
i_2\le \cdots \le i_r\le m-r\}.
\end{equation*} Let $\cal{B}_{m,0} =\{1\}$. 
The set \begin{equation*}\cal{B}_m =\bigcup_{r=0}^m\cal{B}_{m,r}\end{equation*} 
spans $R_m/J_m$.\end{prop}

We prove Proposition \ref{basis} by induction on $m$. The case $m=1$ is trivial, 
but for the inductive step, we need the following lemmas.

Set $J_m =J_{m,0}$ and, for $0<j<m$, define ideals $J_{m,j}$ in $R_m$ 
inductively by
\begin{equation*} J_{m,j} =J_{m,j-1}+\sum_{r=1}^j R_m(Z_1(u)^r)_{m-j} 
=J_m+\sum_{s\ge m-j}\sum_{ 1\le r\le 
m-s}R_m\left(Z_1(u)^r\right)_s.\end{equation*}

\begin{lem}\label{6.5}  If $j\ge 1$, there is a unique homomorphism of 
algebras
\begin{equation*} 
R_{m-j}/J_{m-j}\to R_m/J_{m,j-1}
\end{equation*}
such that $z_i\mapsto z_{i+1}$ for $0\le i<m-j$.
\end{lem}
\begin{pf}  To establish the lemma, we must prove that
\begin{equation}\label{extra}\left(\left(\sum_{i=1}^{m-j}z_iu^i\right)^r\right)_
s\in J_{m,j-1},\end{equation}
for all $r\ge 1$ and $s\ge m-j+1$.

We proceed by induction on $j$. For $j=1$, we must show that
\begin{equation*}\label{j} \left(Z_1(u)^r\right)_s\in J_m \ \ \forall \
 \ s\ge m.\end{equation*}
 If $r=1$, this is trivial from the 
definition of $Z_1(u)$. 
Assume the result for smaller values of $r$. Writing
\begin{equation*} 
Z_0(u) =u(z_0+Z_1(u)),
\end{equation*}
we have 
\begin{equation*} 
\left(Z_0(u)^r\right)_s 
=\left(\sum_{t=0}^r\left(\begin{matrix}r\\t\end{matrix}\right)
z_0^t(Z_1(u)^{r-t})_{s-r}\right). 
\end{equation*}
Take $s=n+r$ with $n\ge m$. Then, $\left(Z_0(u)^r\right)_{n+r}\in J_m$ by the 
definition of $J_m$, so
\begin{equation*}
\sum_{t=0}^r\left(\begin{matrix}r\\t\end{matrix}\right)z_0^t(Z_
1(u)^{r-t})_{n}\in J_m.
\end{equation*}
But, if $t>0$, then $(Z_1(u)^{r-t})_{n}\in J_m$ for all $n\ge m$ by the 
induction 
hypothesis on $r$, so  $(Z_1(u)^{r})_{n}\in J_m$ for all $n\ge m$, thus  
completing the induction on $r$, and establishing (\ref{extra}) when $j=1$.

So now assume that we know (\ref{extra}) for $j-1$. Write
\begin{equation*} 
\sum_{i=1}^{m-j+1}z_iu^i =z_{m-j+1}u^{m-j+1} 
+\sum_{i=1}^{m-j}z_iu^i.
\end{equation*}
By the induction hypothesis on $j$,
\begin{equation*}\left(\left(\sum_{i=1}^{m-j+1}z_iu^i\right)^r\right)_s\in 
J_{m,j-2}\end{equation*}
for all $r\ge 1$, $s\ge m-j+2$. 
Since  $z_{m-j+1}\in J_{m,j-1}$, we can conclude,  by using the binomial 
expansion,  that 
\begin{equation*}
\left(\left(\sum_{i=1}^{m-j} z_iu^i\right)^r\right)_s\in J_{m,j-1},
\end{equation*}
if $r\ge 1$, $s> m-j+1$. Thus, it suffices to prove that
\begin{equation}{\label {5}}
\left(\left(\sum_{i=1}^{m-j}z_iu^i\right)^r\right)_{m-j+1}\in J_{m,j-1}
\end{equation}
for all $r\ge 1$. If $r\le j-1$, we have $(Z_1(u)^r)_{m-j+1}\in J_{m,j-1}$ by 
definition. Further, the elements $z_{m-j+1},\cdots ,z_{m-1}\in J_{m,j-1}$. 
Thus, writing 
\begin{equation*} 
Z_1(u)=\sum_{i=1}^{m-j}z_iu^i+\sum_{i=m-j+1}^{m-1}z_iu^i,
\end{equation*}
and using the binomial expansion,  we see that  (\ref{5}) follows. 

If $r> j-1$, then $r+m-j+1>m$, so
\begin{equation*} \left(Z_0(u)^r\right)_{r+m-j+1}\in J_m\subset 
J_{m,j-1}\end{equation*}
by the definition of $J_m$, we have
\begin{equation*} 
\left((z_0+Z_1(u))^r\right)_{m-j+1}\in J_{m,j-1}.
\end{equation*}
Now, $Z_1(u)_{m-j+1} =z_{m-j+1}\in J_{m,j-1}$ by definition, and so using the 
binomial expansion again and an induction on $r$, we conclude that 
\begin{equation*} 
\left(Z_1(u)^r\right)_{m-j+1}\in J_{m,j-1}.
\end{equation*}
But now the proof is completed as in the case $r<j$.\end{pf}

The proof of the following lemma is elementary.
\begin{lem} Let $r\ge k\ge 0$. Then, the matrix 
\begin{equation*}
\left[ \begin{matrix} 
\left(\begin{matrix}k+1\\1\end{matrix}\right)&\left(\begin{matrix}k+1\\2
\end{matrix}\right)&\cdots&\left(\begin{matrix}k+1\\r-k+1\end{matrix}\right)\\
&&&&\\
\left(\begin{matrix}k+2\\1\end{matrix}\right)&\left(\begin{matrix}k+2\\2
\end{matrix}\right)&\cdots&\left(\begin{matrix}k+2\\r-k+1\end{matrix}\right)\\
.&.&\cdots&.\\
.&.&\cdots&.\\.&.&\cdots&.\\
.&.&\cdots&.\\
\left(\begin{matrix}r+1\\1\end{matrix}\right)&\left(\begin{matrix}r+1\\2
\end{matrix}\right)&\cdots&\left(\begin{matrix}r+1\\r-k+1\end{matrix}\right)
\end{matrix}\right]
\end{equation*}
has determinant $\left(\begin{matrix}r+1\\k\end{matrix}\right)$ and hence is 
invertible.\hfill\qedsymbol
\end{lem}

\begin{lem}{\label{7.9}} For $r\ge t>0$, $0\le s
\le r$, the element $z_0^{r-s}(Z_1(u)^{s+1})_{m-t}$ belongs to the span of
\begin{equation*} 
\{z_0^{t-j}(Z_1(u)^{r-t+j+1})_{m-t}: 1\le j\le t\}.
\end{equation*}
\end{lem}
\begin{pf} We assume that $m>1$, otherwise there is nothing to prove. We 
consider the following equations in $R_m/J_m$:
\begin{equation*} 
z_0^{r-j}(Z_0(u)^{j+1})_{m+j+1-t} =0,\ \ 0<t\le j,
\end{equation*}
i.e.,
\begin{equation*} 
z_0^{r-j}((z_0+Z_1(u))^{j+1})_{m-t} =0,\ \ 0<t\le j,
\end{equation*}
i.e.,
\begin{equation*}
\sum_{i=0}^j\left(\begin{matrix} 
j+1\\i\end{matrix}\right)z_0^{r-j+i}(Z_1(u)^{j+1-i})_{m-t} =0.
\end{equation*}
We must show that these equations, for $j=t, t+1,\cdots ,r$, can be solved for 
the elements $z_0^s(Z_1(u)^{r+1-s})_{m-t}$ with $t\le s\le r$ in terms of those 
with $s<t$. But this follows from the preceding lemma.
\end{pf}

\noindent{\it {Proof of Proposition \ref{basis}.}} The proposition is trivially 
true if $m=1$. Assume now that we know the result for $m-1$. 

For $0\le k<r\le m$, set
\begin{equation*}
\cal{B}_{m,r,k} =\{z_{i_1}z_{i_2}\cdots z_{i_r}\in\cal{B}_{m,r} 
: i_{k+1}\ge 1\}.
\end{equation*}
The proposition obviously follows from

\vskip 6pt

\noindent {\it{Claim}} Let $r\ge k\ge 0$, and let $g\in R_m$ be a homogenous
polynomial of degree $r-k$ in $z_1,z_2,\cdots, z_{m-1}$. Then, $z_0^kg$ is in 
the span of $\cal{B}_{m,r,k}$ modulo $J_m$. 

\medskip

We proceed by induction on $k$. If $k=0$, then by Lemma \ref{6.5} we have a 
homomorphism $R_{m-1}/J_{m-1}\to R_m/J_m$ which sends $z_i\to z_{i+1}$. Clearly, 
$g$ is in the image of this homomorphism and the induction hypothesis on $m$ 
implies that $g\in\cal{B}_{m,r,0}$.

Assume the result for $k-1$. Write
\begin{equation*}
 g = g_0+g_1z_{m-k} + g_2z_{m-k+1}+\cdots +g_{k}z_{m-1},\end{equation*}
where for $0\le j\le k$,  $g_j$ is a polynomial in $z_1,z_2,\cdots 
,z_{m-k+j-1}$.
Now, for $j\ge 0$, we see by Lemma \ref{7.9} that  the element  $z_0^kz_{m-k+j}$ 
is in the span of the sets $\cal{B}_{m,k+1, s}$ with $s<k$. Thus, the element 
$z_0^kz_{m-k+j}g_{j+1}$ can be written as a sum $\sum_{s<k}z_0^sh_{s
j}$, where the $h_{sj}$ are polynomials in $z_1,\cdots ,z_{m-1}$. Hence, by the 
induction on $k$, 
\begin{equation*} z_0^kz_{m-k+j}g_{j+1}\in \cal{B}_{m,r,k},\end{equation*}
for $j\ge 0$. Finally observe that by Lemma \ref{6.5}, $g_0$ is in the image of 
the map $R_{m-k-1}/J_{m-k+j-1}\to R_m/J_{m}$ and hence, by using the induction 
on $m$, we get that
\begin{equation*} g_0\in {\text{span}}(\cal{B}_{m, r-k, 0}) \mod 
J_{m,k+1}.\end{equation*}
Thus,  $z_0^kg_0$ is in the span of $\cal{B}_{m,r,k}$ provided that 
$z_0^kJ_{m,k}$ is also in the span of $\cal{B}_{m,r,k}$, i.e if 
$z_0^k(Z_1(u)^r)_{m-k-1}$ is in the span of $\cal{B}_{m,r,k}$ for all  $s\ge 
m-k$, $1\le r\le m-s$. 
But this follows from Lemma \ref{7.9} again, and the proof of the proposition is 
complete.\hfill\qedsymbol

\end{document}